\documentclass[11pt,reqno]{amsart}
\usepackage{pb-diagram} 
\newtheorem{definition}{Definition}
\newtheorem{theorem}{Theorem}
\newtheorem{proposition}[theorem]{Proposition}
\newtheorem{corollary}[theorem]{Corollary}
\newtheorem{remark}{Remark}
\newtheorem{lemma}[theorem]{Lemma}
\newfont{\bb}{msbm10 at 12pt}

\def\Si{{\Sigma}}     
\def\ga{{\gamma}}     
     
\def\Ga{{\Gamma}}

\def\phi{{\varphi}}

\def\TM{{\mathrm{T}\mathrm{M}}}  
   
\def\M{\mathrm{M}}
\def\T*M{\mathrm{T}^*\mathrm{M}}

\def\H{\mathrm{H}}
\def\R{\mathrm{R}}
\def\gst{{g_{\rm st}}}

\def\L{\mathrm{L}}

\def\L{\mathrm{L}}

\def\C{\mathrm{C}}
\def\D{\mathrm{D}}
\def\hs{\mathbb{S}^n_+}

\def\V{\mathrm{V}}
\def\G{\mathrm{G}}
\def\ov{\overline}
\def\F{\mathrm{F}}
\def\U{\mathrm{U}}
\let\pa\partial     
\let\na\nabla     
     
\DeclareMathAlphabet{\doba}{U}{msb}{m}{n}

\def\T{\mathrm{T}}
\def\Vol{{\mathop{\rm Vol}}}          
\def\Ric{{\mathop{\rm Ric}}}

\def\Spin{{\mathop{\rm Spin}}}     
\def\SO{{\mathop{\rm SO}}}

\def\ker{{\mathop{\rm Ker}}}

\def\Vol{{\mathop{\rm Vol}}}      
      
\def\Ric{{\mathop{\rm Ric}}}

\def\Spin{{\mathop{\rm Spin}}}      
\def\SO{{\mathop{\rm SO}}}

\def\ker{{\mathop{\rm ker}}}

\def\bhlp{\lambda_{\min}^+(\M,[g],\sigma)}

\def\spec{\mathrm{Spec}}
\def\Ker{\mathrm{Ker}}
\def\Coker{\mathcal{C}}
\let\<\langle     
\let\>\rangle

\def\Id{{\mathop{\rm Id}}}

\def\hpm{\mathcal{H}^{\pm}}    
\def\wit{\widetilde}

\let\<\langle     
\let\>\rangle

\long\def\komment#1{} 

\begin{document}


\title[On a spin conformal invariant]{On a spin conformal invariant on manifolds with boundary}   
\author{Simon Raulot}    
\address{Institut {\'E}lie Cartan\\
Universit{\'e} Henri Poincar{\'e}, Nancy I\\
B.P. 239\\
 54506 Vand\oe uvre-L{\`e}s-Nancy Cedex, France}
\email{raulot@iecn.u-nancy.fr}
\date{\today}


\maketitle

\begin{abstract}
Let $\M$ be an $n$-dimensional connected compact manifold with non-empty boundary equipped with a Riemannian metric $g$, a spin structure $\sigma$ and a chirality operator $\Ga$. We define and study some properties of a spin conformal invariant given by:
\begin{eqnarray*}
\lambda_{\mathrm{min}}(\M,\pa\M):=\underset{\ov{g}\in[g]}{\inf}|\lambda_1^\pm(\ov{g})|\Vol(\M,\ov{g})^{\frac{1}{n}},
\end{eqnarray*}

\noindent where $\lambda_1^\pm(\ov{g})$ is the smallest eigenvalue of the Dirac operator under the chiral bag boundary condition $\mathbb{B}^\pm_{\ov{g}}$. More precisely, we show that if $n\geq 2$ then:
\begin{eqnarray*}
\lambda_{\mathrm{min}}(\M,\pa\M)\leq\lambda_{\mathrm{min}}(\hs,\pa\hs).
\end{eqnarray*} 

\end{abstract}


\section{Introduction}


On a compact Riemannian spin manifold without boundary $(\M^n,g)$, several results (see \cite{hijazi:86}, \cite{lott}, \cite{baer:92} or \cite{habilbernd} for instance) have been devoted to a spin conformal invariant independently introduced by Hijazi and Lott and defined in terms of the smallest positive eigenvalue of the Dirac operator. This invariant is given by:
\begin{equation}\label{bhl}
\begin{array}{ccc}
\bhlp & := & \underset{\ov{g}\in[g]}{\inf}\lambda_1^+(\ov{g})\Vol(\M,\ov{g})^{\frac{1}{n}},
\end{array}
\end{equation}  

where $[g]$ and $\sigma$ are respectively the conformal class of the Riemannian metric $g$ and a spin structures on $\M$ (they also define this invariant from the largest negative eigenvalue of the Dirac operator). In \cite{lott} and \cite{ammann}, the authors proved that $\bhlp>0$ using pseudo-differential operators and Sobolev embedding theorems. Moreover, in \cite{ammann} and \cite{amm3} it is shown that:
\begin{eqnarray}\label{large}
\bhlp\leq\lambda_{\min}(\mathbb{S}^n,[\gst],\sigma_{\mathrm{st}})=\frac{n}{2}\omega_n^{\frac{1}{n}},
\end{eqnarray}

where $\mathbb{S}^n$ is the $n$-dimensional Euclidean sphere with standard conformal and spin structure and $\omega_n$ stands for its volume. In \cite{amm4}, Ammann, Humbert and Morel proved that Inequality~(\ref{large}) is strict under some geometrical assumptions. More precisely, they introduced the mass endomorphism on locally conformally flat manifolds which is the constant term of the Green function of the Dirac operator and they showed that if this endomorphism is not identically zero, $n\not\equiv 3\,\text{mod}\,4$ and $\ker(\D)=\{0\}$ then Inequality~(\ref{large}) is strict. The strict inequality has several applications; first using the Hijazi inequality (see \cite{hijazi:86} and \cite{hijazi:91}), it gives a spinorial proof of the Yamabe problem (see \cite{lee.parker:87}); secondly one can obtain solutions of a nonlinear partial differential equation in a way analogous to the Yamabe problem. Note that this equation is critical in the sense that the Sobolev embeddings involved are critical. For a complete review of these results, see \cite{ahm}.

In this paper, we define and study an analogous spin conformal invariant in the context of compact manifolds with boundary. On such manifolds, the Dirac operator has an infinite dimensional kernel and a closed image with finite codimension so we have to impose some boundary conditions on the restriction to the boundary $\pa\M$ of the spinor fields on $\M$. More details on elliptic boundary conditions for the Dirac operator can be found in \cite{bw:93}, \cite{hijazi.montiel.roldan:01} or \cite{fs}. In order to define a well-posed spin conformal invariant, we have to choose some adapted boundary condition, i.e. a conformally invariant boundary condition. The chiral bag boundary condition is such a condition and so on a compact Riemannian spin manifold with boundary we let:
\begin{eqnarray*}
\lambda_{\mathrm{min}}(\M,\pa\M):=\underset{\ov{g}\in[g]}{\inf}|\lambda_1^\pm(\ov{g})|\Vol(\M,\ov{g})^{\frac{1}{n}},
\end{eqnarray*}

where $\lambda_1^\pm(\ov{g})$ is the smallest eigenvalue of the Dirac operator under the chiral bag boundary condition $\mathbb{B}^\pm_{\ov{g}}$ (see Section~\ref{cbbc}), $g$ (resp. $[g]$) is a Riemannian metric (resp. the conformal class of $g$) and $\sigma$ a spin structure on $\M$. The main result of this article is the following:
\begin{theorem}\label{main}
Let $(\M^n,g,\sigma)$ be a connected compact Riemannian spin manifold with smooth boundary equipped with a chirality operator $\Ga$. Then:
\begin{eqnarray}\label{largeb}
\lambda_{\mathrm{min}}(\M,\pa\M)\leq\lambda_{\mathrm{min}}(\hs,\pa\hs)=\frac{n}{2}\Big(\frac{\omega_n}{2}\Big)^{\frac{1}{n}},
\end{eqnarray}

where $\omega_n$ stands for the volume element of the standard sphere $\mathbb{S}^n$. 
\end{theorem}

The proof of this theorem is given in four steps. First in Section~\ref{vc}, we give a variational characterization of this invariant. Then in Section~\ref{ts}, we will compute explicitly its value on the standard hemisphere and in particular we will construct Killing spinor fields which satisfy the chiral bag boundary condition. Section~\ref{spt} is devoted to the construction of a trivialization of the spinor bundle over the manifold $\M$ equipped with a Riemannian metric which sends (locally) a spinor field on the flat space to a spinor field over an open set of $\M$. Finally from a Killing spinor field and this trivialization, we will obtain an adapted test spinor to evaluate in the variational characterization of $\lambda_{\min}(\M,\pa\M)$. In \cite{sr4}, we give a sufficient condition for a certain type of manifolds for which Inequality~(\ref{largeb}) is strict. This condition is based on the construction of the Green function for the Dirac operator under the chiral bag boundary condition.\\

This invariant is closely related to the Yamabe invariant involved in the Yamabe problem on manifolds with boundary. This problem could be stated as follows: given $(\M^n,g)$ a compact Riemannian manifold of dimension $n\geq3$, does there exist a metric $\ov{g}$ conformal to $g$ for which the scalar curvature is constant and the mean curvature is zero? This problem has been solved in many cases by Escobar \cite{escobar:92}. Solving the Yamabe problem is equivalent to finding a positive smooth function $f$ of the following boundary value problem:
$$\left\lbrace
\begin{array}{ll}
\L_g f:=\frac{4(n-1)}{n-2}\Delta_g f+\R_g f= Cf^{N-1} & \quad\text{on}\;\;\M\\
\mathrm{B}_g f_{|\pa\M}:=\frac{2}{n-2}\frac{\pa f}{\pa\nu}+\H_g f=0 & \quad\text{along}\;\;\pa\M
\end{array}
\right.$$

where $\R_g$ (resp. $\H_g$) is the scalar (resp. mean) curvature of $(\M,g)$ (resp. $(\pa\M,g_{|\pa\M})$), $N=\frac{2n}{n-2}$ and $C$ is a constant. The existence of a smooth positive solution to this system is based on the study of the Yamabe invariant defined by 
\begin{eqnarray*}
\mu(\M,\pa\M):= \underset{u\in\mathrm{C}^1(\M),u\neq 0}{\mathrm{inf}} \frac{\int_{\M}\big(\frac{2}{n-2}|\na u|^2+\frac{1}{2(n-1)}\R_g u^2\big)dv(g)+\int_{\pa\M}\H_g u^2
  ds(g)}{\Big(\int_{\M}|u|^N dv(g)\Big)^{\frac{2}{N}}}.
\end{eqnarray*}
 
This number only depends on the conformal class of $g$ and if $\mu(\M,\pa\M)\geq 0$ it can be expressed as
\begin{eqnarray*}
\mu(\M,\pa\M)=\inf\mu_{1}(\L_{\ov{g}})\Vol(\M,\ov{g})^{\frac{2}{n}},
\end{eqnarray*}
 
where the infimum is taken over all metrics $\ov{g}$ conformal to $g$ and where $\mu_1(\L_g)$ is the first eigenvalue of the eigenvalue problem:
\begin{equation}\label{pvpy}
\left\lbrace
\begin{array}{ll}
\L_g u= \mu_1(\L_g)u & \quad\text{on}\;\;\M\\
\mathrm{B}_g u_{|\pa\M}=0 & \quad\text{along}\;\;\pa\M
\end{array}
\right.
\end{equation}

In \cite{escobar:92}, Escobar shows that:
\begin{eqnarray}\label{pablo}
\mu(\M,\pa\M)\leq\mu(\hs,\pa\hs)=n(n-1)\Big(\frac{\omega_n}{2}\Big)^{\frac{2}{n}}.
\end{eqnarray}

Moreover, he proves that if this inequality is strict then the Yamabe problem is solved. To point out the relation between $\lambda_{\mathrm{min}}(\M,\pa\M)$ and this problem, let us recall some results. If $n\geq 3$, then the Hijazi inequality on manifolds with boundary \cite{sr} relates $\lambda_1^\pm(g)$ to the first eigenvalue $\mu_1(\L_g)$:
\begin{eqnarray*}
\lambda_1^\pm(g)^2\geq\frac{n}{4(n-1)}\mu_{1}(\L_g)
\end{eqnarray*}

and equality holds if and if only $(\M,g)$ is isometric to a half-sphere. Using the H\"older inequality gives:
\begin{eqnarray}\label{hijbord}
\lambda_{\mathrm{min}}(\M,\pa\M)^2\geq\frac{n}{4(n-1)}\mu(\M,\pa\M).
\end{eqnarray}

{\bf Acknowledgements:} I would like to thank Oussama Hijazi and Emmanuel Humbert for their support. I am also very grateful to Marc Herzlich and Sebasti\'an Montiel for their remarks and their suggestions.


\section{Chiral bag boundary condition}\label{cbbc}


In this section, we give some standard facts about compact Riemannian spin manifold with boundary and chiral bag boundary condition. This condition has been introduced in \cite{hawking} to prove some positive mass theorems for black holes (see also \cite{herzlich1}). Note that this condition has also been studied in a serie of papers by Gilkey, Kirsten and others (see \cite{gilkey}, \cite{gilkey1} and \cite{gilkey2}). For more details on boundary conditions for the Dirac operator, we refer to \cite{bw:93} or \cite{hijazi.montiel.roldan:01} .

Let $(\M^n,g)$ be a compact Riemannian spin manifold with boundary, denote by $\Sigma_g(\M)$ its spinor bundle, $\na$ its Levi-Civita connection and ``$\cdot$'' its Clifford multiplication. The spinor bundle is endowed with a natural Hermitian scalar product denoted by $\<\,,\,\>$ compatible with $\na$ and with the Clifford multiplication. The Dirac operator $\D_g$ is then the first order elliptic differential operator acting on $\Sigma_g(\M)$ locally given by:
\begin{eqnarray*}
\D_g\varphi=\sum_{i=1}^n e_i\cdot\na_{e_i}\varphi,
\end{eqnarray*}

for all $\varphi\in\Ga\big(\Sigma_g(\M)\big)$ and where $\{e_1,...,e_n\}$ is a local $g$-orthonormal frame of the tangent bundle. Since the boundary is an oriented hypersurface of $\M$, there exists a unit vector field $\nu$ normal to $\pa\M$ and then the boundary is itself a spin manifold. We denote by $\mathbf{S}_g:=\Sigma_g(\M)_{|\pa\M}$ the restriction of the spinor bundle of $\M$ to the boundary. This bundle is also endowed with a Levi-Civita connection $\na^{\mathbf{S}}$, a Clifford multiplication and a natural compatible Hermitian scalar product.

Now suppose that there exists on the manifold $\M$ a chirality operator $\Ga$, i.e. an endomorphism of the spinor bundle satisfying:
\begin{equation}\label{poc}
\begin{array}{rc}
\Ga^2=\Id, & \<\Ga\varphi,\Ga\psi\>=\<\varphi,\psi\>\\
\na_X(\Ga\psi)=\Ga(\na_X\psi), & X\cdot\Ga\psi=-\Ga(X\cdot\psi)
\end{array}
\end{equation}
 
for all $\varphi$, $\psi\in\Ga\big(\Si_g(\M)\big)$ and $X\in\Ga(\TM)$. Since the endomorphism $\nu\cdot\Ga$ is involutive, we have an eigenbundle decomposition $\mathbf{S}_g=\V^+\oplus\V^-$ where $\V^{\pm}$ is the eigensubbundle associated with the eigenvalue $\pm 1$. We can then check that the projection:
$$\begin{array}{lccl}
\mathbb{B}^{\pm}_g: & \L^2(\mathbf{S}_g) & \longrightarrow & \L^2(\V^{\pm})\\
 & \varphi & \longmapsto & \frac{1}{2}(\Id\pm \nu\cdot\Ga)\varphi,
\end{array}$$ 

defines an elliptic boundary condition for the Dirac operator called the chiral bag boundary condition. Indeed, in \cite{fs} the authors show that:
\begin{eqnarray}
\D_g:\hpm_g=\{\varphi\in\H_1^2\;/ \;\mathbb{B}^{\pm}_g(\varphi_{|\pa\M})=0\}\longrightarrow\L^2\big(\Si_g(\M)\big)
\end{eqnarray}

defines a self-adjoint Fredholm operator, i.e. the kernel of $\D_g$ is finite dimensional and its image is closed in $\L^2\big(\Si_g(\M)\big)$. Then we have an $\L^2$-orthogonal splitting:
\begin{eqnarray}\label{decort}
\mathcal{H}^{\pm}_g=\Ker^{\pm}(\D_g)\oplus\Coker^{\pm}_g,
\end{eqnarray}

where $\Coker^{\pm}_g$ is the $\L^2$-orthogonal of $\Ker^{\pm}(\D_g)$ in $\mathcal{H}^{\pm}_g$. It is easy to check that $\Coker^{\pm}_g$ is closed in $\mathcal{H}^{\pm}_g$. Moreover, under this boundary condition, the spectrum of the Dirac operator $\D_g$ with domain $\mathcal{H}^{\pm}_g$ consists of entirely isolated real eigenvalues with finite multiplicity and it admits a spectral resolution:
\begin{eqnarray}
\L^{2}\big(\Si_g(\M)\big)=\underset{\lambda^{\pm}(g)\in\spec^{\pm}(\D_g)}{\bigoplus}\mathcal{N}_{\lambda^{\pm}(g)}(\D_g),
\end{eqnarray}

where $\spec^{\pm}(\D_g)$ is the spectrum under the boundary condition $\mathbb{B}^{\pm}_g$ and $\mathcal{N}_{\lambda^{\pm}(g)}(\D_g)$ is the eigenspace associated with the eigenvalue $\lambda^{\pm}(g)$. In \cite{sr}, we show that this set of eigenvalues satisfies the Hijazi inequality~(\ref{hijbord}). We are now ready to define the chiral bag invariant.
\begin{definition}
The chiral bag invariant is given by:
\begin{eqnarray*}
\lambda_{\min}^{\pm}(\M,\pa\M):=\underset{\ov{g}\in[g]}{\inf}|\lambda^{\pm}_1(\ov{g})|\Vol(\M,\ov{g})^{\frac{1}{n}},
\end{eqnarray*}

where $\lambda_1^{\pm}(\ov{g})$ is the smallest eigenvalue of the Dirac operator under the chiral bag boundary condition $\mathbb{B}^{\pm}_{\ov{g}}$.
\end{definition}

\begin{remark}
{\rm This definition seems to depend on the boundary condition chosen $\mathbb{B}^{+}_{g}$ or $\mathbb{B}^{-}_{g}$, however it doesn't. Indeed, we are going to check that $\spec^{+}(\D_g)=-\spec^{-}(\D_g)$ and then $\lambda_{\min}^{+}(\M,\pa\M)=\lambda_{\min}^{-}(\M,\pa\M)$. Let $\varphi\in\Gamma\big(\Si_g(\M)\big)$ be an eigenspinor for the Dirac operator under the $\mathbb{B}^{-}_{g}$ boundary condition, then it satisfies the eigenvalue problem:
$$\left\lbrace
\begin{array}{ll}
\D_g\varphi = \lambda\varphi & \qquad\text{on}\;\M\\
\nu\cdot\Ga\varphi_{|\pa\M}=\varphi_{|\pa\M} & \qquad\text{along}\;\pa\M
\end{array}
\right.$$

We can then decompose the eigenspinor following the chirality decomposition, i.e. $\varphi=\varphi^++\varphi^-$ where $\varphi^{\pm}=\frac{1}{2}(1\pm\Ga)\varphi$. Since the Dirac operator sends positive spinors to negative ones and conversely, an easy calculation leads to $\D_g\widetilde{\varphi}=-\lambda\widetilde{\varphi}$, where $\widetilde{\varphi}=\varphi^+-\varphi^-$. Using the boundary condition, the chirality decomposition and since Clifford multiplication by $\nu$ interchanges the chirality of spinor fields, we can easily check that $\nu\cdot\varphi^{\pm}_{|\pa\M}=\pm\varphi^{\mp}_{|\pa\M}$. We deduce that:
\begin{eqnarray*}
\nu\cdot\Ga\widetilde{\varphi}_{|\pa\M}=\nu\cdot\varphi^{+}_{|\pa\M}+\nu\cdot\varphi^{-}_{|\pa\M}=\varphi^-_{|\pa\M}-\varphi^+_{|\pa\M}=-\widetilde{\varphi}_{|\pa\M},
\end{eqnarray*}
 
and then $-\lambda\in\spec^+(\D_g)$. In particular, we have $\Ker^{+}(\D_g)\simeq\Ker^{-}(\D_g)$ and Clifford multiplication by $\nu$ sends $\Coker^{+}_g$ on $\Coker^{-}_g$.}
\end{remark}

Taking this into account, we will denote by $\lambda_{\min}(\M,\pa\M)$ the chiral bag invariant in the rest of this paper and we will consider the $\mathbb{B}^-_g$ condition.


\section{Variational characterization}\label{vc}


The aim of this section is to give a variational characterisation of the chiral bag invariant. More precisely, we prove that:
\begin{proposition}\label{vcbhl}
We have:
\begin{eqnarray}
\lambda_{\min}(\M,\pa\M):=\underset{\ov{g}\in[g]}{\inf}|\lambda_1^-(\ov{g})|\Vol(\M,\ov{g})^{\frac{1}{n}}=
\underset{\varphi\in\Coker^-_g}{\inf}\Big\{\frac{\big(\int_{\M}|\D_g\varphi|^{\frac{2n}{n+1}}dv(g)\big)^{\frac{n+1}{n}}}{\big|\int_{\M}\mathrm{Re}\<\D_g\varphi,\varphi\>dv(g)\big|}\Big\}.
\end{eqnarray}

\end{proposition}

In order to prove this proposition, we need a variational characterization of the first eigenvalue of the Dirac operator under the chiral bag boundary condition. So we first show:
\begin{proposition}\label{vcpv}
The square of the first eigenvalue of the Dirac operator is given by:
\begin{eqnarray*}
\lambda_1^-(g)^2=\inf_{\varphi\in\Coker^-_g}\Big\{\frac{\int_{\M}|\D_g\varphi|^2 dv(g)}{\int_{\M}|\varphi|^2 dv(g)}\Big\}.
\end{eqnarray*}
\end{proposition}

We omit the proof of this proposition which can be seen as a direct application of the Rayleigh quotient. However, we refer to \cite{mathese} for a rigorous proof. We can easily deduce from the Cauchy-Schwarz inequality and from the preceding proposition the following result which is very useful for the variational characterization of $\lambda_{\mathrm{min}}(\M,\pa\M)$.
\begin{corollary}\label{pvd}
The first eigenvalue of the chiral bag boundary condition is given by:
\begin{eqnarray*}
|\lambda_1^-(g)|=\underset{\psi\in\Coker^-_g}{\inf}\Big\{\frac{\int_{\M}|\D_g\psi|^2dv(g)}{\big|\int_{\M}\mathrm{Re}\<\D_g\psi,\psi\>dv(g)\big|}\Big\}.
\end{eqnarray*}
\end{corollary}

From now on, all quantities which depend of the metric are written with their reference metric. We now briefly recall some conformal aspects of spin geometry. Let $\ov{g}=f^2 g$ be a metric in the conformal class of $g$, then there exists an isomorphism $\F$ between their respective spinor bundles $\Si_g(\M)$ and $\Si_{\ov{g}}(\M)$. We can also relate the Dirac operators $\D_g$ and $\D_{\ov{g}}$ (see \cite{hitchin:74} or \cite{hijazi:86}) by the formula:
\begin{eqnarray}\label{covdir}
\D_{\ov{g}}(\F(\psi))=f^{-\frac{n+1}{2}}\F(\D_g(f^{\frac{n-1}{2}}\psi)).
\end{eqnarray}

The chiral bag boundary condition transforms nicely under conformal change of metrics. Indeed, if there exists a chirality operator acting on $\Sigma_g(\M)$, then the map
\begin{eqnarray*}
\ov{\Ga}:=\F\circ\Ga\circ\F^{-1}:\Sigma_{\ov{g}}(\M)\longrightarrow\Sigma_{\ov{g}}(\M)
\end{eqnarray*}

defines a chirality operator on the spinor bundle over $\M$ endowed with the metric $\ov{g}$.Thus the orthogonal projection
$$\begin{array}{lccl}
\mathbb{B}^{\pm}_{\ov{g}}: & \L^2(\mathbf{S}_{\ov{g}}) & \longrightarrow & \L^2(\ov{\V}^{\pm})\\
 & \varphi & \longmapsto & \frac{1}{2}(\Id\pm \ov{\nu}\;\ov{\cdot}\ov{\Ga})\varphi,
\end{array}$$ 

also defines elliptic boundary condition for the Dirac operator. Then using the conformal covariance of the Dirac operator and that of the chiral bag boundary condition, we easily check that:
\begin{eqnarray*}
\psi\in\Ker(\D_g)\quad\Longleftrightarrow\quad f^{-\frac{n-1}{2}}\F(\psi)\in\Ker(\D_{\ov{g}}).
\end{eqnarray*}

Its $\L^2$-orthogonal complement also transforms naturally as:
\begin{eqnarray*}
\psi\in\Coker^-_g\quad\Longleftrightarrow\quad f^{-\frac{n+1}{2}}\F(\psi)\in\Coker^-_{\ov{g}}.
\end{eqnarray*}

{\it Proof of Proposition~\ref{vcbhl}:} The proof is closely related to the one given in \cite{ammann}. First we note that since the Dirac operator under the chiral bag boundary condition admits a self-adjoint $\L^2$-extention, we can identify the orthogonal of the kernel with the image of $\D_{\ov{g}}$. With a slight abuse in the notations, let $\Psi=f^{-\frac{n+1}{2}}\F(\varphi)\in\Coker^-_{\ov{g}}$ and then using Corollary~\ref{pvd} and the preceding discussion, we can write:
\begin{eqnarray}\label{class}
|\lambda_1^-(\ov{g})| & = & \underset{\Psi\in\Coker^-_{\ov{g}}}{\inf}\Big\{\frac{\int_{\M}|\Psi|^2dv(\ov{g})}{\big|\int_{\M}\mathrm{Re}\<\Psi,\D^{-1}_{\ov{g}}\Psi\>dv(\ov{g})\big|}\Big\}\nonumber\\ \nonumber \\ 
& = &  \underset{\Psi\in\Coker^-_{g}}{\inf}\Big\{\frac{\int_{\M}|f^{-\frac{n+1}{2}}\F(\varphi)|^2 f^n dv(g)}{\big|\int_{\M}\mathrm{Re}\<f^{-\frac{n+1}{2}}\F(\varphi),\D^{-1}_{\ov{g}}\big(f^{-\frac{n+1}{2}}\F(\varphi)\big)\> f^n dv(g)\big|}\Big\}\nonumber\\  \nonumber\\ 
& = & \underset{\varphi\in\Coker^-_g}{\inf}\Big\{\frac{\int_{\M} f^{-1}|\varphi|^2 dv(g)}{\big|\int_{\M}\mathrm{Re}\<\varphi,\D^{-1}_{g}\varphi\big)\> dv(g)\big|}\Big\}.
\end{eqnarray}

Now note that, using the H\"older inequality, we have:
\begin{eqnarray*}
\Big(\int_{\M}|\varphi|^{\frac{2n}{n+1}}dv(g)\Big)^{\frac{n+1}{n}} \leq \Big(\int_{\M}f^{-1}|\varphi|^2dv(g)\Big)^{\frac{n}{n+1}}\mathrm{vol}(\M,\ov{g})^{\frac{1}{n}},
\end{eqnarray*}

for all $f\in\C^{\infty}_{+}(\M)$, $\varphi\in\Coker^-_g$ and with equality if and only if $f=c\,|\varphi|^{\frac{2}{n+1}}$ or $\varphi\equiv 0$. We then obtain that:
\begin{eqnarray*}
\underset{f\in\C^{\infty}_{+}(\M)}{\inf}\Big\{\big(\int_{\M}f^{-1}|\varphi|^2dv(g)\big)^{\frac{n}{n+1}}\mathrm{vol}(\M,\ov{g})^{\frac{1}{n}}\Big\} & \geq &  ||\varphi||^{2}_{\L^{2n/(n+1)}},
\end{eqnarray*}

for $\varphi\in\Coker^-_g$. However, if the spinor field $\varphi$ has no zeros we can let $f=|\varphi|^{\frac{2}{n+1}}$ and then we have equality in the preceding inequality. Otherwise we can easily find a sequence $f_k:\M\rightarrow\mathbb{R}^{*}_{+}$ such that if $\ov{g}_k=f^2_kg$ then:
\begin{eqnarray*}
\big(\int_{\M}f_{k}^{-1}|\varphi|^2dv(g_1)\big)^{\frac{n}{n+1}}\mathrm{vol}(\M,\ov{g}_{k})^{\frac{1}{n}}\underset{k\rightarrow\infty}{\longrightarrow}  ||\varphi||^{2}_{\L^{2n/(n+1)}},\quad\text{for all}\;\varphi\in\Coker^-_g.
\end{eqnarray*}

Using all these arguments and the fact that $\ov{g}=f^2 g$, the chiral bag invariant is then given by:
\begin{eqnarray*}
\lambda_{\min}(\M,\pa\M) & = & \underset{f\in\C^{\infty}_{+}(\M)}{\inf}|\lambda_1^-(\ov{g})|\mathrm{vol}(\M,\ov{g})^{\frac{1}{n}}\\ \\
& = & \underset{f\in\C^{\infty}_{+}(\M)}{\inf}\;\underset{\varphi\in\Coker^-_g}{\inf}\Big\{\frac{\big(\int_{\M}f^{-1}|\varphi|^2dv(g)\big)}{\big|\int_{\M}\mathrm{Re}\<\varphi,\D^{-1}_g\varphi\>dv(g)\big|}\mathrm{vol}(\M,\ov{g})^{\frac{1}{n}}\Big\} \\ \\
& = & \underset{\varphi\in\Coker^-_g}{\inf}\;\underset{f\in\C^{\infty}_{+}(\M)}{\inf}\Big\{\frac{\big(\int_{\M}f^{-1}|\varphi|^2dv(g)\big)}{\big|\int_{\M}\mathrm{Re}\<\varphi,\D^{-1}_g\varphi\>dv(g)\big|}\mathrm{vol}(\M,\ov{g})^{\frac{1}{n}}\Big\}\\ \\
& = & \underset{\varphi\in\Coker^-_g}{\inf}\Big\{\frac{\big(\int_{\M}|\varphi|^{\frac{2n}{n+1}}dv(g)\big)^{\frac{n+1}{n}}}{\big|\int_{\M}\mathrm{Re}\<\varphi,\D^{-1}_g\varphi\>dv(g)\big|}\Big\}.
\end{eqnarray*}
\hfill$\square$


\section{The case of the hemisphere}\label{ts}


In this section, we consider the particular case where the manifold is the hemisphere endowed with its standard spin structure and conformal class. In fact, we show that on $\hs$, we can construct Killing spinor fields satisfying  the chiral bag boundary condition. Thus we will be able to compute explicitly the value of $\lambda_{\min}(\hs,\pa\hs)$. First we prove the following result:
\begin{proposition}\label{testspinor}
On the half-space $\mathbb{R}^n_+$ endowed with the standard Euclidian metric $\xi$, there exists a spinor field $\psi^{\pm}\in\Ga\big(\Si_{\xi}(\mathbb{R}^n_+)\big)$ which satisfies the following boundary problem:
$$\left\lbrace
\begin{array}{ll}
\D_\xi\psi^{\pm}= \pm\frac{n}{2} f\psi^{\pm} & \rm{on}\;\mathbb{R}^n_+\\\
\mathbb{B}^-_\xi(\psi^\pm_{|\pa\mathbb{R}^n_+})=0 \quad & \rm{along}\;\pa\mathbb{R}^n_+
\end{array}
\right.$$

where $\D_\xi$ is the Dirac operator on $\mathbb{R}^n_+$ and  $f$ is the real-valued function given by $f(x)=\frac{2}{1+r^2}$ with $r^2=x_1^2+...+x_n^2$ if $x=(x_1,...,x_n)\in\mathbb{R}^n_+$. Moreover, we have the following relations:
\begin{eqnarray}\label{killingspinors}
|\D_\xi\psi^{\pm}|^2=f^{n+1}\quad\text{and}\quad|\psi^{\pm}|^2=f^{n-1}.
\end{eqnarray}
\end{proposition}

{\it Proof:}
Fix $\Phi_0$ a parallel spinor on $\Si_{\xi}(\mathbb{R}^n_+)$ and for $x=(x_1,...,x_n)\in\mathbb{R}^n_+$, we define:
\begin{eqnarray*}
\psi^{\pm}(x)=\frac{1}{\sqrt{2}}f^{\frac{n}{2}}(x)(1\mp x)\cdot\Phi_0(x).
\end{eqnarray*} 

An easy computation using the parallelism of $\Phi_0$ and the relation $\partial_i\cdot \partial_j+\partial_j\cdot \partial_i=-2\delta_{ij}$ for all $1\leq i,j\leq n$ leads to:
\begin{eqnarray*}\label{de}
\D_\xi\psi^+ (x)= -\frac{n}{2\sqrt{2}}f^{\frac{n}{2}+1}(x)\,x\cdot\Phi_0+\frac{n}{2\sqrt{2}}f^{\frac{n}{2}}(x)\big(2-f(x) r^2\big)\Phi_0.
\end{eqnarray*}

Note that:
\begin{eqnarray*}
2-f(x) r^2=f(x)
\end{eqnarray*}

and then:
\begin{eqnarray}\label{de}
\D_\xi\psi^+ (x)=\frac{n}{2} f(x)\psi^+(x).
\end{eqnarray}
 
The preceding calculation depends only on the parallelism of the spinor field $\Phi_0$, so we now show that this spinor can be chosen satisfying the chiral bag boundary condition. In fact, if we let $\Psi_0=\frac{1}{2}(\Phi_0+\widetilde{\nu}\cdot\Ga\Phi_0)$ where $\widetilde{\nu}$ is a smooth extension of the normal field $\nu$, then, since $\pa\mathbb{R}^n_+$ is totally geodesic in $\mathbb{R}^n_+$, the spinor field $\Psi_0$ is also parallel hence it also satisfies (\ref{de}). Now using the properties~(\ref{poc}) of the chirality operator, we have:
\begin{eqnarray*}
\nu\cdot\Ga\psi^+_{|\pa\mathbb{R}^n_+}=\frac{1}{2}f^{\frac{n}{2}}(x)\nu\cdot\Ga\big((1- x)\cdot(\Phi_0+\nu\cdot\Ga\Phi_0)\big)_{|\pa\mathbb{R}^n_+}=\psi^+_{|\pa\mathbb{R}^n_+}.
\end{eqnarray*}

Finally, we can easily compute that this spinor field satisfies the relations (\ref{killingspinors}). 
\hfill$\square$\\

Using this result, we can now construct a Killing spinor field on the standard hemisphere $\hs$ satisfying the chiral bag boundary condition. Indeed, we have:
\begin{corollary}\label{spinchi}
On the standard hemisphere $\hs$, there exists a Killing spinor $\psi^{\pm}\in\Ga\big(\Si_{\gst}(\hs)\big)$ satisfying the following eigenvalue problem:
\begin{equation}\label{sk1}
\left\lbrace
\begin{array}{ll}
\D^{\hs}(\varphi^{\pm})=\pm\frac{n}{2}\varphi^{\pm}\quad & \rm{on}\;\hs\\
\mathbb{B}^-_{\gst}(\varphi^{\pm}_{|\pa\mathbb{S}^n_+})=0 \quad & \rm{along}\;\pa\hs,
\end{array}
\right.
\end{equation}

where $\D^{\hs}$ is the Dirac operator on $\hs$.
\end{corollary}

{\it Proof:} First recall that if $q\in\pa\hs$, then the stereographic projection of pole $q$ gives an isometry between $\big(\mathbb{R}^n_+,f^2\,g_{{\rm eucl}}\big)$ and $\big(\hs\setminus \{q\},\gst\big)$, where $f$ is the function given in Proposition~\ref{testspinor}. The bundle isomorphism $\F$ described in Section~\ref{vc} allows to identify the spinor bundle over $\mathbb{R}^n_+$ with that over $\hs\setminus\{q\}$. So if we let $\varphi^{\pm}:=f^{-\frac{n-1}{2}}\F(\psi^{\pm})$ where $\psi^\pm$ is the spinor field constructed in Proposition~\ref{testspinor}, the conformal covariance of the Dirac operator~(\ref{covdir}) gives:
\begin{eqnarray*}
\D^{\hs}(\varphi^{\pm}) = \pm\frac{n}{2}\varphi^{\pm}
\end{eqnarray*}

and the one of the chiral bag boundary condition leads to $\mathbb{B}^-_{\gst}(\varphi^\pm_{|\pa\hs})=0$ on $\hs\setminus \{q\}$. Now we show that this spinor field extends to the whole hemisphere. For $\varepsilon>0$ let $\eta_\varepsilon:\U\rightarrow [0,1]$ be a cut-off function with support in an open set $\U$ of $\hs$ and with $q\in\U$. More precisely, we suppose that $\eta_\varepsilon\equiv 1$ on $\mathrm{B}^+_\varepsilon(q)$, $\textrm{supp}\,(\eta_\varepsilon)\subset\mathrm{B}^+_{2\varepsilon}(q)$ and $|\nabla\eta_\varepsilon|\leq C / \varepsilon$ where $C$ is a positive real number. Observe now that for $\Phi\in\mathcal{H}^-_{\gst}$ we have:
\begin{eqnarray}\label{rst}
\int_{\mathrm{U}}\<\varphi^\pm,\D^{\hs}\Phi\>dv(g_{\rm{st}}) & = & \int_{\mathrm{U}}\<\varphi^\pm,\D^{\hs}\big(\eta_\varepsilon\Phi+(1-\eta_\varepsilon)\Phi\big)\>dv(g_{\rm{st}})\nonumber\\
 & = & \int_{\mathrm{U}}\<\varphi^\pm,\D^{\hs}\big((1-\eta_\varepsilon)\Phi\big)\>dv(g_{\rm{st}})+\int_{\mathrm{U}}\<\varphi^\pm,\eta_\varepsilon\D^{\hs}\Phi\>dv(g_{\rm{st}})\nonumber\\ 
 & & +\int_{\mathrm{U}}\<\varphi^\pm,\nabla\eta_\varepsilon\cdot\Phi\>dv(g_{\rm{st}}).
\end{eqnarray}

Since $\varphi^\pm$ satisfies the chiral bag boundary condition on $\hs\setminus\{q\}$, we get:
\begin{eqnarray*}
\int_{\mathrm{U}}\<\varphi^\pm,\D^{\hs}\big((1-\eta_\varepsilon)\Phi\big)\>dv(g_{\rm{st}}) =  \pm\frac{n}{2}\int_{\mathrm{U}}\<\varphi^\pm,(1-\eta_\varepsilon)\Phi\>dv(g_{\rm{st}}) \underset{\varepsilon\rightarrow 0}{\longrightarrow} \pm\frac{n}{2}\int_{\U}\<\varphi^\pm,\Phi\>dv(g_{\rm{st}}).
\end{eqnarray*}

On the other hand, an estimation of the second term in Identity~(\ref{rst}) leads to:
\begin{eqnarray*}
\big|\int_\U\<\varphi^\pm,\eta_\varepsilon\D^{\hs}\Phi\>dv(g_{\rm{st}})\big|\leq ||\varphi^\pm||_{\L^2(\mathrm{B}^+_{2\varepsilon})}||\D^{\hs}\Phi||_{\L^2(\mathrm{B}^+_{2\varepsilon})}.
\end{eqnarray*}

It is then easy to conclude that the righthand side of this inequality goes to $0$ when $\varepsilon\rightarrow 0$. Let's now look at the last term in (\ref{rst}). We have:
\begin{eqnarray*}
\big|\int_{\mathrm{U}}\<\varphi^\pm,\nabla\eta_\varepsilon\cdot\Phi\>dv(g_{\rm{st}})\big|\leq \mathrm{vol}\big(\mathrm{B}^+_{2\varepsilon}(q)\big)^{\frac{1}{2}}\,||\na\eta_{\varepsilon}||_{\mathrm{C^0}}\,||\varphi^\pm||_{\mathrm{L}^2}\,||\Phi||_{\mathrm{L}^2}\leq C||\varphi^\pm||_{\L^2(\mathrm{B}^+_{2\varepsilon})}\,\varepsilon^{\frac{n}{2}-1}
\end{eqnarray*}

and since $n\geq 2$, it converges to $0$ when $\varepsilon\rightarrow 0$. Hence, we have proved that the spinor field $\varphi^\pm$ satisfies, in a weak sense on $\U$, the eigenvalue boundary problem (\ref{sk1}). The classical regularity theorems allow to conclude that it is satisfied, in a strong sense, on $\hs$.
\hfill$\square$\\

We can thus easily compute the value of the chiral bag invariant on the hemisphere. Indeed we have the following result:
\begin{corollary}
The chiral bag invariant on the hemisphere $\mathbb{S}^n_+$ is given by:
\begin{eqnarray}
\lambda_{\min}(\hs,\pa\hs)=\frac{n}{2}\Big(\frac{\omega_{n}}{2}\Big)^{\frac{1}{n}}.
\end{eqnarray}
\end{corollary}

{\it Proof:} First, we note that using the Hijazi inequality~(\ref{hijbord}) we have:
\begin{eqnarray*}
\lambda_{\min}(\hs,\pa\hs)^2\geq\frac{n}{4(n-1)}\mu(\mathbb{S}^n_+)=\frac{n^2}{4}\Big(\frac{\omega_{n}}{2}\Big)^{\frac{2}{n}}.
\end{eqnarray*}

Moreover (see Corollary~\ref{spinchi}), since a Killing spinor is an eigenspinor for the Dirac operator under the chiral bag boundary condition, we obtain:
\begin{eqnarray*}
\lambda_{\min}(\hs,\pa\hs):=\underset{\ov{g}\in[\gst]}{\inf}|\lambda_1(\ov{g})|\Vol(\mathbb{S}^n_+,\ov{g})^{\frac{1}{n}}\leq|\lambda_1(\gst)|\Vol(\mathbb{S}^n_+,\gst)^{\frac{1}{n}}=\frac{n}{2}\Big(\frac{\omega_{n}}{2}\Big)^{\frac{1}{n}}
\end{eqnarray*}

and we can easily conclude.
\hfill $\square$


\section{A trivialization of the spinor bundle}\label{spt}


In order to prove Theorem~\ref{main}, we have to construct a trivialization of the spinor bundle over $\M$ endowed with a Riemannian metric around a boundary point. This trivialization arises from a bundle isomorphism introduced by Bourguignon and Gauduchon \cite{bourguignon} to identify spinors on a Riemannian spin manifold endowed with two distinct metrics and from an adapted chart of the manifold around a boundary point, the Fermi coordinates. We follow more particularly \cite{amm3}.\\

Let $q$ be a boundary point and let $(x_1,...,x_{n-1})$ be normal coordinates on $\pa\M$ at this point. Let $t\mapsto\ga(t)$ be the geodesic leaving from $(x_1,...,x_{n-1})$ in the orthogonal direction to $\pa\M$ and parametrized by arc length. Then:
$$\begin{array}{lclc}
\mathcal{F}_q : & \U\subset\T_{q}\M & \longrightarrow & \V\subset\M \\
 & (x_1,...,x_{n-1},t) & \longmapsto & m
\end{array}$$
 
are called the Fermi coordinates at $q\in\pa\M$. Moreover in these coordinates, the arc length $ds^2$ is given by:
\begin{eqnarray}\label{mcf}
ds^2=dt^2+g_{ij}(x,t)dx^idx^j,
\end{eqnarray}

for all $1\leq i,j\leq n-1$ and $(x,t)\in \U$. Now let:
$$\begin{array}{lcll}
\mathrm{G}: & \V\subset\M & \longrightarrow & \mathrm{S}^{2}_{+}(n,\mathbb{R})=\{\rm{A}\in\mathcal{M}_{n}(\mathbb{R})\,/\,\rm{A}\text{ symmetric and positive-definite}\}\\
 & m & \longmapsto & \mathrm{G}_{m}=\big(g_{ij}(m)\big)_{1\leq i,j\leq n}
\end{array}$$

the smooth map which associates to any point $m\in\V$ the matrix of the coefficients of the metric $g$ at this point in the basis $\{\pa_1,...,\pa_{n-1},\pa_t\}$. In Fermi coordinate, we can also write:
$$\mathrm{G}_m=\left(
\begin{array}{cc}\widetilde{\mathrm{G}}_m & 0 \\
0 & 1
\end{array}
\right),\quad\text{with}\;\widetilde{\mathrm{G}}_{m}=\big(g_{ij}(m)\big)_{1\leq i,j\leq n-1}.$$

Since the manifold $\M$ is Riemannian, the metric is positive-definite and symmetric at each point $m\in\V$, hence $\mathrm{G}_m$ too. In a same way, the matrix $\widetilde{\mathrm{G}}_m$ is also positive-definite and symmetric. So there exists $\widetilde{\mathrm{B}}_m\in\mathrm{S}^{2}(n-1,\mathbb{R})$ such that:
\begin{eqnarray*}
\widetilde{\mathrm{B}}_m^2=\widetilde{\mathrm{G}}_{m}^{-1},
\end{eqnarray*} 

which depends smoothly on $m$. If we let:
$$\mathrm{B}_m=\left(
\begin{array}{cc}
\widetilde{\mathrm{B}}_m & 0 \\
0 & 1
\end{array}
\right)\in\mathrm{S}^{2}(n,\mathbb{R}),$$

then this matrix satisfies $\mathrm{B}_m^2=\mathrm{G}_{m}^{-1}$. Note that for all $X,Y\in\mathbb{R}^{n}_{+}$, we have:
\begin{eqnarray}\label{iso}
^{t}\big(\mathrm{B}_m X\big)\mathrm{G}_m\big(\mathrm{B}_m Y\big)=\;^{t}X\;\mathrm{I}_n Y=\xi(X,Y),
\end{eqnarray}

where $\xi$ stands for the standard euclidian metric on $\mathbb{R}^{n}_{+}$. We then have an isomorphism:
$$\begin{array}{llll}
\mathrm{B}_m: & \big( \T_{\mathcal{F}^{-1}_{q}(m)}\U\simeq\mathbb{R}^{n}_{+},\xi \big) & \longrightarrow & \big(\T_m\V,g_m\big),
\end{array}$$

which, by construction, depends smoothly on $m$. We can now identify the $\SO_n$-principal bundles $\SO(\U,\xi)$ and $\SO(\V,g)$ of oriented $\xi$ and $g$-orthonormal frames of $(\U,\xi)$ and $(\V,g)$. In fact, the following diagram is commutative:
$$\begin{diagram}
\node{\SO(\U,\xi)}\arrow{e,t}{\zeta}\arrow{s,l}{}\node{\SO(\V,g)}\arrow{s,r}{}\\
\node{\U\subset\T_{q}\M}\arrow{e,b}{\mathcal{F}_{q}}\node{\V\subset\M}
\end{diagram}$$

where $\zeta$ is given by the natural action of $\mathrm{B}$ on $\SO(\U,\xi)$. This diagram commutes with the right action of $\SO_n$, then the map $\zeta$ can be lifted to:
$$\begin{diagram}
\node{\Spin(\U,\xi)}\arrow{e,t}{\widetilde{\zeta}}\arrow{s,l}{}\node{\Spin(\V,g)}\arrow{s,r}{}\\
\node{\U\subset\T_{q}\M}\arrow{e,b}{\mathcal{F}_{q}}\node{\V\subset\M}
\end{diagram}$$

Hence, we obtain an identification between the spinor bundles $\Sigma_{\xi}(\U)$ and $\Sigma_g(\V)$ given by:
\begin{equation}
\begin{array}{clc}\label{iden}
\Sigma_{\xi}(\U):=\Spin(\U,\xi)\times_{\rho_n}\Sigma_{n} & \longrightarrow & \Sigma_g(\V):=\Spin(\V,g)\times_{\rho_n}\Sigma_{n} \\
 \psi=[s,\varphi]& \longmapsto & \overline{\psi}=[\widetilde{\zeta}(s),\varphi].
\end{array}
\end{equation}

where $(\rho_n,\Si_n)$ is the complex spinor representation. In the same way, we can identify the boundary spinor bundles $\mathbf{S}_{\xi}(\U):=\Sigma_{\xi}(\U)_{|\U\cap\pa\mathbb{R}^n_+}$ and $\mathbf{S}_{g}(\V):=\Sigma_{g}(\V)_{|\V\cap\pa\M}$. \\

Now we are going to relate the Dirac operator $\D_g$ acting on sections of $\Sigma_g(\V)$ with the Dirac operator $\D_\xi$ acting on those of $\Si_{\xi}(\U)$. First, let:
\begin{eqnarray*}
e_i:=b^{j}_{i}\pa_j,
\end{eqnarray*}

where $b^{j}_{i}$ are the coefficients of the matrix $\mathrm{B}_m$ and then $\{e_1,...,e_n\}$ is a local orthonormal frame of $(\T\V,g)$. We can also suppose that the unit vector field $e_n$ is the inner unit vector field $\nu:=\pa_t$ normal to $\pa\M$ and then we have $b^{j}_{n}=\delta^{j}_{n}$ and $b^{n}_{i}=\delta^{n}_{i}$ for all $1\leq i,j\leq n$. Denote by $\na$ (resp. $\ov{\na}$) the Riemannian and spinorial Levi-Civita connection on $(\U,\xi)$ (resp. $(\V,g)$) and the Christoffel of the second kind $\wit{\Ga}_{ij}^k$ are defined by:
\begin{eqnarray*}
\wit{\Ga}_{ij}^k=g(\ov{\na}_{e_i} e_j,e_k).
\end{eqnarray*} 

We can then easily check that Clifford multiplications on $\Sigma_{\xi}(\U)$ and $\Sigma_g(\V)$ are related by:
\begin{eqnarray*}
\ov{\pa_i\cdot\psi}=e_i\cdot\ov{\psi},
\end{eqnarray*}

for all $\psi\in\Sigma_{\xi}(\U)$.  Moreover,  for a spinor field $\psi\in\Ga\big(\Si_{\xi}(\U)\big)$ then $\ov{\psi}\in\Ga\big(\Si_g(\V)\big)$ and by construction of the spinorial Levi-Civita connection, we have:
\begin{eqnarray}\label{levicivita}
\ov{\na}_{e_i}\overline{\psi}=e_i(\psi)+\frac{1}{4}\sum_{1\leq j,k\leq n}\wit{\Gamma}_{ij}^{k}\,e_j\cdot e_k\cdot\overline{\psi}.
\end{eqnarray}

Using the local expression of the Dirac operator $\D_g$, we get:
\begin{eqnarray}\label{identdir}
\D_g\ov{\psi}=\ov{\D_\xi\psi}+\sum_{i,j=1}^{n}\big(b_i^j-\delta_i^j\big)\ov{\pa_i\cdot\na_{\pa_j}\psi}+\frac{1}{4}\sum_{1 \leq i,j,k\leq n}\wit{\Ga}_{ij}^{k}\,e_i\cdot e_j\cdot e_k\cdot\ov{\psi}.
\end{eqnarray}

We can then prove the following statement:
\begin{proposition}\label{devdir}
If $\D_\xi$ and $\D_g$ are the Dirac operators acting respectively on $\Si_{\xi}(\U)$ and $\Si_g(\V)$, then we have
\begin{eqnarray}
\D_g\overline{\psi}=\overline{\D_\xi\psi}+\sum_{i,j=1}^{n}\big(b_i^j-\delta_i^j\big)\ov{\pa_i\cdot\na_{\pa_j}\psi}+\mathrm{W}\cdot\overline{\psi}+\mathrm{T}\cdot\overline{\psi}+\widetilde{\nu}\cdot\mathrm{Z}\cdot\overline{\psi}-\frac{n-1}{2}\mathrm{H}_t\,\widetilde{\nu}\cdot\overline{\psi},
\end{eqnarray}

where $\widetilde{\nu}\in{\Gamma(\T\M)}$ is a local extention of the inner normal vector field $\nu\in{\Gamma(\T\M_{|\pa\M})}$ and where $\mathrm{W}\in\Gamma(\Lambda^3(\T^*\V))$, $\mathrm{T}\in\Gamma(\T^*\V)$ and  $\mathrm{Z}\in\Gamma(\Lambda^2(\T^*\V))$ are given by
\begin{eqnarray*}
\mathrm{W} & = & \frac{1}{4}\sum_{\underset{i\neq j\neq k}{1\leq i,j,k\leq n-1}}b^{r}_{i}\pa_{r}(b^{l}_{j})(b^{-1})^{k}_{l}\,e_i\cdot e_j \cdot e_k\\
\mathrm{T} & = &\frac{1}{4}\sum_{1\leq i,j\leq n-1}(\widetilde{\Gamma}_{ij}^{i}-\widetilde{\Gamma}_{ii}^{j})\,e_j\\
\mathrm{Z} & = &\frac{1}{4}\sum_{\underset{i\neq j}{1\leq i,j\leq n-1}}\big(\pa_n(b^{l}_{i})(b^{-1})^{j}_{l}+b^{r}_{j}\Gamma_{rn}^{l}(b^{-1})^{i}_{l}\big)\,e_i\cdot e_j,
\end{eqnarray*}

with $\H_t$ the mean curvature of $\pa\M_t:=\{\Theta_q(x,t)\;/ \;t\;\rm{is\;constant}\}$ and where, for any point $m\in\V$, the coefficients $(b^{-1})^{k}_{l}$ are the coefficients of the inverse matrix of $\mathrm{B}_m$.
\end{proposition}

{\it Proof:} 
In order to compute the second term of formula (\ref{identdir}), we decompose the sum into tangential and normal parts. An easy calculation using the fact that $\wit{\Ga}^n_{nn}=0$ leads to
\begin{eqnarray*}
\sum_{1\leq i,j,k\leq n}\wit{\Ga}_{ij}^{k}\,e_i\cdot e_j\cdot e_k\cdot\ov{\psi} & = &\sum_{i,j,k=1}^{n-1}\wit{\Ga}_{ij}^{k}\,e_i\cdot e_j\cdot e_k\cdot\ov{\psi}+\wit{\nu}\cdot\sum_{i,j=1}^{n-1}\wit{\Ga}_{ni}^{j}\,e_i\cdot e_j \cdot\ov{\psi}\\
& & -\wit{\nu}\cdot \sum_{i,j=1}^{n-1}\wit{\Ga}_{in}^{j}\,e_i\cdot e_j\cdot\ov{\psi}+\wit{\nu}\cdot\sum_{i,j=1}^{n-1}\wit{\Ga}_{ij}^{n}\,e_i\cdot e_j\cdot\ov{\psi}\\
& & - \sum_{i=1}^{n-1}\big(\wit{\Ga}_{nn}^{i}-\wit{\Ga}_{ni}^{n}+\wit{\Ga}_{in}^{n}\big)e_i\cdot\ov{\psi}.
\end{eqnarray*}

Now note that we have $\wit{\Ga}^{k}_{ij}=-\wit{\Ga}^{j}_{ik}$, $\wit{\Ga}^{i}_{nn}=0$ by construction and :
\begin{eqnarray*}
\wit{\Ga}_{in}^{n}=e_i(g(\widetilde{\nu},\widetilde{\nu}))-g(e_i,\ov{\na}_{\widetilde{\nu}}\widetilde{\nu})=-\wit{\Ga}_{in}^{n}
\end{eqnarray*}

hence $\wit{\Ga}_{in}^{n}=0$ for all $1\leq i\leq n$. Then the preceding equality gives:
\begin{eqnarray*}
\sum_{i,j,k=1}^{n}\wit{\Ga}_{ij}^{k}\,e_i\cdot e_j\cdot e_k\cdot\ov{\psi} = 
\sum_{i,j,k=1}^{n-1}\wit{\Ga}_{ij}^{k}\,e_i\cdot e_j\cdot e_k\cdot\ov{\psi}+\wit{\nu}\cdot\sum_{i,j=1}^{n-1}\big(\wit{\Ga}_{ni}^{j}-2\wit{\Ga}_{in}^{j}\big)e_i\cdot e_j\cdot\ov{\psi}.
\end{eqnarray*}

However, the last term of this expression can be simplified using the fact that $\wit{\Ga}_{ni}^{i}=-\wit{\Ga}_{ni}^{i}=0$ to give:
\begin{eqnarray*}
\wit{\nu}\cdot\sum_{\underset{i\neq j}{i,j=1}}^{n-1}\big(\wit{\Ga}_{ni}^{j}+\wit{\Ga}_{jn}^{i}-\wit{\Ga}_{in}^{j}\big)\,e_i\cdot e_j\cdot\ov{\psi}-2(n-1)\H_t\,\ov{\psi},
\end{eqnarray*}

where $\H_t=\frac{1}{n-1}\sum_{i=1}^{n-1}g(-\ov{\na}_{e_i}\wit{\nu},e_i)$ is the mean curvature of $\pa\M_t$. A direct computation gives:
\begin{eqnarray*}
\sum_{i,j,k=1}^{n-1}\wit{\Ga}_{ij}^{k}\,e_i\cdot e_j\cdot e_k\cdot\ov{\psi} & = & \sum_{i\not{=} j\neq k\neq i}\wit{\Ga}_{ij}^{k}\,e_i\cdot e_j\cdot e_k\cdot\ov{\psi}+\sum_{i,j=1}^{n-1}\big(\wit{\Ga}_{ij}^{i}-\wit{\Ga}_{ii}^{j}\big)\,e_j\cdot\ov{\psi},
\end{eqnarray*}

then combining all the above results, we have:
\begin{eqnarray*}
\D_g\psi & = & \ov{\D_\xi\psi}+\sum_{i,j=1}^{n}\big(b_i^j-\delta_i^j\big)\ov{\pa_i\cdot\na_{\pa_j}\psi}+\overbrace{\frac{1}{4}\sum_{\underset{i\neq j\neq k}{i,j,k=1}}^{n-1}\wit{\Ga}_{ij}^{k}\,e_i\cdot e_j\cdot e_k\cdot\ov{\psi}}^{(1)}+\frac{1}{4}\sum_{i,j=1}^{n-1}\big(\wit{\Ga}_{ij}^{i}-\wit{\Ga}_{ii}^{j}\big)\,e_j\cdot\ov{\psi}\nonumber\\
& & -\frac{n-1}{2}\mathrm{H}_t\,\widetilde{\nu}\cdot\ov{\psi}\nonumber+\underbrace{\frac{1}{4}\,\widetilde{\nu}\cdot\sum_{\underset{i\neq j}{i,j=1}}^{n-1}\big(\wit{\Ga}_{ni}^{j}+\wit{\Ga}_{jn}^{i}-\wit{\Ga}_{in}^{j}\big)\,e_i\cdot e_j\cdot\ov{\psi}}_{(2)}.
\end{eqnarray*}

We are now going to give the expansion of the preceding expression in terms of the coefficients $(b_i^j)_{1\leq i,j\leq n}$. First consider $(1)$; by construction of the $b^j_i$, we have $\wit{\Ga}_{ij}^k e_k=\wit{\Ga}_{ij}^k b^l_k\pa_l$ and otherwise:
\begin{eqnarray*}
\wit{\Ga}_{ij}^{k} e_k=\ov{\na}_{e_i}e_j=b^r_i\ov{\na}_{\pa_r}(b^s_j\pa_s)=b^r_i\pa_r(b^s_j)\pa_s+b^r_i b^s_j\Ga_{rs}^{l}\pa_l,
\end{eqnarray*}

where the Christoffel symbols $\Ga_{rs}^{l}$ are given by $\Ga_{rs}^{l}=g(\ov{\na}_{\pa_r}\pa_s,\pa_l)$. We then have:
\begin{eqnarray*}
\wit{\Gamma}_{ij}^{k} e_k=\big(b^r_i\pa_r(b^l_j)+b^r_i b^s_j\Gamma_{rs}^{l}\big)\pa_l,
\end{eqnarray*}

and so:
\begin{eqnarray}\label{christoffel}
\wit{\Ga}_{ij}^{k}=\big(b^r_i\pa_r(b^l_j)+b^r_i b^s_j\Ga_{rs}^{l}\big)(b^{-1})^k_l. 
\end{eqnarray}

The term $(1)$ is hence given by:
\begin{eqnarray*}
\sum_{\underset{i\neq j\neq k}{i,j,k=1}}^{n-1}\wit{\Ga}_{ij}^{k}\, e_i\cdot e_j\cdot e_k\cdot\ov{\psi} & = & \sum_{\underset{i\neq j\neq k}{i,j,k=1}}^{n-1}\big(b^r_i\pa_r(b^l_j)+b^r_i b^s_j\Ga_{rs}^{l}\big)(b^{-1})^k_l\,e_i\cdot e_j\cdot e_k\cdot\ov{\psi}. 
\end{eqnarray*}

Using the symmetry of the Christoffel symbols and the relation $e_i\cdot e_j=-e_j\cdot e_i$ for $i\neq j$, we easily check that $(1)$ is given by:
\begin{eqnarray*}
\frac{1}{4}\sum_{\underset{i\neq j\neq k}{i,j,k=1}}^{n-1}\wit{\Ga}_{ij}^{k}\,e_i\cdot e_j\cdot e_k\cdot\ov{\psi}=\frac{1}{4}\sum_{\underset{i\neq j\neq k}{1\leq i,j,k\leq n-1}}b^{r}_{i}\pa_{r}(b^{l}_{j})(b^{-1})^{k}_{l}\,e_i\cdot e_j\cdot e_k\cdot\ov{\psi}=\mathrm{W}\cdot\ov{\psi}
\end{eqnarray*}

with $\mathrm{W}\in\Gamma(\Lambda^3(\T^*\V))$. For $(2)$, note that we have $\wit{\Ga}_{ni}^{j}=\big(\pa_n(b_i^l)+b^r_i\Ga_{nr}^l\pa_l\big)(b^{-1})^j_l$ and:
\begin{eqnarray*}
\wit{\Ga}_{jn}^i e_i=\wit{\Ga}_{jn}^i b^l_i\pa_l=\ov{\na}_{e_j}\widetilde{\nu}=b^r_j\ov{\na}_{\pa_r}\pa_t=b^r_j\Ga_{rn}^l\pa_l,
\end{eqnarray*}

so identifiying the two parts leads $\wit{\Ga}_{jn}^i=b^r_j\Ga_{rn}^l(b^{-1})^i_l$. The symmetry of the Christoffel symbols allows to conclude:
\begin{eqnarray*}
(2)  & = & \frac{1}{4}\widetilde{\nu}\cdot\sum_{\underset{i\neq j}{i,j=1}}^{n-1}\big(\pa_t(b^l_i)(b^{-1})^j_l+b^r_j\Ga_{rn}^l(b^{-1})^i_l\big)\;e_i\cdot e_j\cdot\ov{\psi}=\widetilde{\nu}\cdot\mathrm{Z}\cdot\ov{\psi},
\end{eqnarray*}

where $\mathrm{Z}\in\Ga\big(\Lambda^2(\T^*\V)\big)$. 
\hfill$\square$


\section{Manifolds of dimension $n\geq 3$}\label{mwnub}


In this section, we use all the preceding results to prove the following theorem:
\begin{theorem}\label{main2}
Let $(\M,g,\sigma)$ be a connected compact Riemannian spin manifold of dimension $n\geq 3$ with non-empty smooth boundary 
equipped with a chirality operator $\Ga$. Then
\begin{eqnarray}\label{largeb1}
\lambda_{\min}(\M,\pa\M)\leq\lambda_{\min}(\hs,\pa\hs)=\frac{n}{2}\Big(\frac{\omega_n}{2}\Big)^{\frac{1}{n}},
\end{eqnarray}

where $\omega_n$ stands for the volume element of the standard sphere $\mathbb{S}^n$. 
\end{theorem}
 
The proof of this theorem is based on the computation of the first terms of the Taylor development of the Dirac operator in the trivialization introduced in Section~\ref{spt}. The final step consists of the construction of a test spinor satisfying the chiral bag boundary condition using the one of Section~\ref{ts}.


\subsection{Expansion of the metric}\label{dm1}


Here we give the development of the metric in Fermi coordinates around a boundary point in order to obtain the Taylor expansion of the Dirac operator in the trivialization of Section~\ref{spt} in terms of curvatures. Let $q\in\pa\M$ and for any point $m$ with $\V\ni m\simeq(x,t)\in\U$, denote by $r=|(x,t)|$ the distance from $q$ to $m$. Recall that in Fermi coordinates in the neighborhood of $q$, we have the following development of the inverse $g^{-1}$ of the metric $g$ (see \cite{esc}):
\begin{equation}
\begin{array}{lll}\label{dm}
g^{ij}  & = & \delta^{ij}+2 h^{ij}(q)t-\frac{1}{3} R^{i\quad j}_{\;\alpha\beta}(q)\;x^{\alpha}x^{\beta}+g^{ij}_{\;\, ,t\alpha}(q)\;x^{\alpha}t \\
& & + \big(3h^{im}(q)h_{m}^{\quad j}(q)+\wit{R}^{\,i\quad j}_{\;\; n n}(q)\big)\;t^2+O(r^3),
\end{array}
\end{equation}

for all $1\leq i,j\leq n-1$ and where $h^{ij}=g^{ik}g^{jl}h_{kl}$ ($h_{kl}$ are the components of the second fundamental form), $\wit{R}^{\,i\quad j}_{\;\; \alpha\beta}$ (resp. $R^{i\quad j}_{\;\alpha\beta}$) those of the Riemann curvature tensor of the manifold $\M$ (resp. of the boundary $\pa\M$) and $g^{ij}_{\;\, ,t\alpha}:=\pa_\alpha\pa_t g^{ij}$. Then we write:
\begin{eqnarray*}
\G^{-1}_m=\Id+\G_1+\G_2+\G_3+\G_4+O(r^3),
\end{eqnarray*}

with:
$$\begin{array}{ll}
(\G_1)^{ij}=2 h^{ij}(q)t, & (\G_2)^{ij}= -\frac{1}{3} R^{i\quad j}_{\;\alpha\beta}(q)\;x^{\alpha}x^{\beta}\\ 
(\G_3)^{ij}=g^{ij}_{\;\, ,t\alpha}(q)\;x^{\alpha}t, & (\G_4)^{ij}=\big(3h^{im}(q)h_{m}^{\quad j}(q)+\wit{R}^{\,i\quad j}_{\;\; n n}(q)\big)\;t^2.
\end{array}$$

Using the fact that $\G_m\G^{-1}_m=\mathrm{Id}$, an easy computation gives:
\begin{equation}
\begin{array}{lll}\label{Devmet}
g_{ij} & = & \delta_{ij}-2 h_{ij}(q)t+\frac{1}{3}R_{i\alpha\beta j}(q) x^{\alpha}x^{\beta}+g_{ij,t\alpha}(q)x^{\alpha} t\\
& & +\big(h_{im}(q)h^{m}_{\;\;\; j}(q)-\wit{R}_{injn}(q)\big)t^2+O(r^3).
\end{array}
\end{equation}

In order to have the Taylor development of $\mathrm{W}\in\Gamma(\Lambda^3(\T^*\V))$, $\mathrm{T}\in\Gamma(\T^*\V)$ and  $\mathrm{Z}\in\Gamma(\Lambda^2(\T^*\V))$ of Proposition~\ref{devdir}, we first consider the development of the $b^j_i$-coefficients. Writing:
\begin{eqnarray*}
\mathrm{B}_m=\Id+\mathrm{B}_1+\mathrm{B}_2+\mathrm{B}_3+\mathrm{B}_4+\mathrm{B}_5+O(r^3),
\end{eqnarray*}

with:
$$\begin{array}{ll}
(\mathrm{B}_1)_{ij}=\mathrm{B}_{ij\alpha}x^{\alpha}, & (\mathrm{B}_2)_{ij}=\mathrm{B}_{ijt}t\\ 
(\mathrm{B}_3)_{ij}=\mathrm{B}_{ij\alpha\beta} x^{\alpha}x^{\beta}, & (\mathrm{B}_4)_{ij}=\mathrm{B}_{ij\alpha t}x^{\alpha}t,\\
(\mathrm{B}_5)_{ij}=\mathrm{B}_{ij t^2}t^2, &
\end{array}$$

and the relation $\mathrm{B}_m^2=\G^{-1}_m$ yields to:
\begin{equation}
\begin{array}{lll}\label{DEVBij}
b^{j}_{i} & = & \delta^{j}_{i}+ h^{ij}(q) t-\frac{1}{6} R^{i\quad j}_{\;\alpha\beta}(q)\;x^{\alpha}x^{\beta} + \frac{1}{2} g^{ij}_{\;\, ,t\alpha}(q)\;x^{\alpha}t \\
& & + \frac{1}{2} \big(2h^{im}(q)h_{m}^{\quad j}(q)+\wit{R}^{\,i\quad j}_{\;\; n n}(q)\big)\,t^2+O(r^3).
\end{array}
\end{equation}

An analogous computation gives the following inverse $(b^{-1})^{j}_{i}$-coefficients expansion:
\begin{equation}
\begin{array}{l}\label{DEVB-1ij}
(b^{-1})^{j}_{i}=\delta^{j}_{i}- h^{ij}(q) t+\frac{1}{6} R^{i\quad j}_{\;\alpha\beta}(q)\;x^{\alpha}x^{\beta}- \frac{1}{2} g^{ij}_{\;\, ,t\alpha}(q)\;x^{\alpha}t - \frac{1}{2} \wit{R}^{\,i\quad j}_{\;\; n n}(q)\,t^2+O(r^3)
\end{array}
\end{equation}

We also need the first derivaties developments:
\begin{equation}
\begin{array}{lll}\label{derbij}
\pa_k b^j_i & = & -\frac{1}{6}\big(R^{i\quad j}_{\;k\alpha}(q)+R^{i\quad j}_{\;\alpha k}(q)\big)\;x^{\alpha}+ \frac{1}{2} g^{ij}_{\;\, ,t k}(q)\;t+O(r^2)
\end{array}
\end{equation}

and:
\begin{equation}
\begin{array}{lll}\label{derbij}
\pa_t b^j_i & = & h^{ij}(q) + \frac{1}{2} g^{ij}_{\;\, ,t\alpha}(q)\;x^{\alpha} + \big(2h^{im}(q)h_{m}^{\quad j}(q) + \wit{R}^{\,i\quad j}_{\;\; n n}(q)\big)\,t + O(r^2).
\end{array}
\end{equation}

We are now ready to prove the following statement:
\begin{proposition}\label{devdir1}
The fields $\mathrm{W}\in\Gamma(\Lambda^3(\T^*\V))$, $\mathrm{Z}\in\Gamma(\Lambda^2(\T^*\V))$ and $\mathrm{T}\in\Gamma(\T^*\V)$ given in Proposition~\ref{devdir} satisfy: $|\mathrm{W}|=O(r^2)$, $|\mathrm{Z}|=O(r^2)$ and $|\mathrm{T}|=O(r)$.
\end{proposition}

{\it Proof:}
First consider the $3$-form $\mathrm{W}$; using Proposition~\ref{devdir} and Identity~(\ref{derbij}), and since $\pa_k(b^l_j)$ has no constant term, we can observe that any term of order $1$ in $b^{r}_{i}\pa_{r}(b^{l}_{j})(b^{-1})^{k}_{l}$ is a product of $0$-order term of $b^{r}_{i}$ and of a term of $1$-order of $\pa_{r}(b^{l}_{j})$, then we have:
\begin{eqnarray*}
\mathrm{W} & = & \frac{1}{4}\sum_{\underset{i\neq j\neq k}{1\leq i,j,k\leq n-1}}\big(\pa_{i}(b^{k}_{j})+O(r^2)\Big)\,e_i\cdot e_j\cdot e_k.
\end{eqnarray*}

Moreover since $\pa_{i}(b^{k}_{j})=\pa_{i}(b^{j}_{k})$ and for $j\neq k$, $e_j\cdot e_k=-e_k\cdot e_j$ then 
\begin{eqnarray*}
\sum_{\underset{i\neq j\neq k}{1\leq i,j,k\leq n-1}}\pa_{i}(b^{k}_{j})\,e_i\cdot e_j\cdot e_k=0
\end{eqnarray*}

and $|\mathrm{W}|=O(r^2)$. We now investigate the $2$-form $\mathrm{Z}$; because of the expression of $\mathrm{Z}$, we first develop the Christoffel symbols $\Gamma_{rn}^{l}$ for $1\leq r,l\leq n-1$. Recall that:
\begin{eqnarray*}
\Gamma_{rn}^{l}=\frac{1}{2}g^{lk}\big(\pa_r g_{nk}+\pa_n g_{rk}-\pa_k g_{rn}\big)
\end{eqnarray*}

and since $1\leq k,r\leq n-1$ then $g_{nk}=g_{rn}=0$, so we have $\Gamma_{rn}^{l}=\frac{1}{2}g^{lk}\pa_n g_{rk}$. Using (\ref{Devmet}), we deduce that:
$$\begin{array}{lll}
\pa_n g_{rk} & = & -2 h_{rk}(q)+g_{rk,t\alpha}(q)x^{\alpha}+2\big(h_{rm}(q)h^{m}_{\;\;\; k}(q)-\wit{R}_{rnkn}(q)\big)t+O(r^2).
\end{array}$$

Similarly, we compute:
$$\begin{array}{lll}
\Gamma_{rn}^{l} & = & -h_{rl}(q)+\frac{1}{2}g_{rl,t\alpha}(q)x^{\alpha}+\big(\mathrm{K}_{rl}(q)-2h^{lk}(q) h_{rk}(q) \big)t+O(r^2),
\end{array}$$

where $\mathrm{K}_{rk}=h_{rm}h^{m}_{\;\;\; k}-\wit{R}_{rnkn}$  and finally with the help of (\ref{DEVBij}) and (\ref{DEVB-1ij}), we obtain:
$$\begin{array}{lll}
b^r_j\Gamma_{rn}^{l}(b^{-1})^i_l & = & -h_{ij}(q)+\frac{1}{2}g_{ij,t\alpha}(q)x^{\alpha}+\big(\mathrm{K}_{ji}(q)-2h^{ik}(q) h_{kj}(q) \big)t+O(r^2).
\end{array}$$

The first term $\pa_n(b^l_i)(b^{-1})^j_l$ in $\mathrm{Z}$ is given by:
$$\begin{array}{lll}
\pa_n(b^l_i)(b^{-1})^j_l & = & \big(h^{il}(q)+\frac{1}{2}g^{il}_{\;\; ,t\alpha}(q)x^{\alpha}+\widetilde{\mathrm{K}}^{il}(q)t+O(r^2)\big)\big(\delta^j_l-h^{lj}(q)t+O(r^2)\\ \\ 
& = & h^{ij}(q)+\frac{1}{2}g^{ij}_{\;\; ,t\alpha}(q)x^{\alpha}+\big(\widetilde{\mathrm{K}}^{ij}(q)-h^{il}(q)h^{jl}(q)\big)t+O(r^2).
\end{array}$$  

where $\wit{\mathrm{K}}^{il}=2h^{im}h_{m}^{\quad l}+\wit{R}^{\,i\quad l}_{\;\; n n}$. Combining these two developments leads $|\mathrm{Z}|=O(r^2)$. To conclude we examine the $1$-form $\T$; first we are going to relate the Christoffel symbols $\wit{\Ga}_{ij}^k$ with the Christoffel symbols of first kind $\Gamma^{k}_{ij}=g(\ov{\na}_{\pa_i}\pa_j,\pa_k)$. Using the classical formula:
\begin{eqnarray*}
\Gamma_{ij}^{k}=\frac{1}{2}g^{kl}\big(\pa_i g_{jl}+\pa_j g_{il}-\pa_l g_{ij}\big) 
\end{eqnarray*}

and a simple derivation on the development (\ref{Devmet}) of the metric tensor components gives:
\begin{eqnarray*}
\pa_k g_{ij}=\frac{1}{3}\big(R_{ik\alpha j}(q)+R_{i\alpha k j}(q)\big)x^{\alpha}+g_{ij,tk}(q)t+O(r^2)
\end{eqnarray*}

and then it yields to:
$$\begin{array}{lll}
\Gamma_{ij}^{k} & = &  \frac{1}{6}\Big(R_{ji\alpha k}(q)+R_{j\alpha ik}(q)+R_{ij\alpha k}(q)+R_{i\alpha jk}(q)-R_{ik\alpha j}(q)-R_{i\alpha kj}(q)\Big)x^{\alpha}\\ \\
& & +\frac{1}{2}\big(g_{jk,ti}(q)+g_{ik,tj}(q)-g_{ij,tk}(q)\big)t+O(r^2). 
\end{array}$$

The symmetries of the Riemannian curvature tensor give:
\begin{eqnarray*}
\Gamma_{ij}^{k} = -\frac{1}{3}\Big(R_{ik\alpha j}(q)+R_{i\alpha kj}(q)\Big)x^{\alpha}+\frac{1}{2}\big(g_{jk,ti}(q)+g_{ik,tj}(q)-g_{ij,tk}(q)\big)t+O(r^2),
\end{eqnarray*}

Using Proposition~\ref{devdir} and the fact that $\Gamma_{ij}^{k}$ has no constant term, we have:
\begin{eqnarray*}
\wit{\Ga}_{ij}^{k} & = & \pa_i(b^k_j)+\Gamma_{ij}^{k}+O(r^2).
\end{eqnarray*}

The $1$-form $\mathrm{T}$ then clearly satisfies $\mathrm{T}=O(r)$. Moreover, we can give an explicit computation of the development of $\mathrm{T}$. Indeed, we get:
\begin{eqnarray*}
\T=\sum_{j=1}^{n-1}\Big(-\frac{1}{4}\Ric(q)_{\alpha j}x^{\alpha}-\frac{1}{2}\wit{\Ric}(q)_{tj} t + O(r^2)\Big)\;e_j,
\end{eqnarray*}

where we used the Codazzi equation:
\begin{eqnarray}\label{codazzi}
R_{ijkn}=h_{ik,j}-h_{jk,i}
\end{eqnarray}

and the formula $g_{ij,t}=-2 h_{ij}$.
\hfill$\square$\\


\subsection{The Estimate}\label{e1}


In this section, we give the proof of Inequality~(\ref{largeb1}). However, we need the following lemma which gives the existence of an adapted test spinor in the trivialization constructed in Section \ref{spt}:
\begin{lemma}\label{trivialisationbord}
Let $\U$ and $\V$ be the open sets of the trivialization constructed in Section~\ref{spt} and let $\Phi_0\in\Ga\big(\Sigma_{\xi}(\U)\big)$ be a parallel spinor such that $\nu\cdot\Ga\ov{\Phi}_0(q)=\ov{\Phi}_0(q)$ for one point $q\in\pa\V\cap\M$. Then we have: 
\begin{eqnarray*}
\nu\cdot\Ga\ov{\Phi}_{0|\pa\V\cap\M}=\ov{\Phi}_{0|\pa\V\cap\M},
\end{eqnarray*}

i.e. $\mathbb{B}^-_g(\ov{\Phi}_{0|\pa\V\cap\M})=0$.
\end{lemma}

 {\it Proof:} Consider the function $f(p)=|\nu\cdot\Ga\ov{\Phi}_0-\ov{\Phi}_0|^2(p)$ defined on $\V$, then we have to show that $f$ vanishes along the boundary $\pa\V\cap\M$. However, for $1\leq i\leq n-1$ we have:
 \begin{eqnarray*}
 e_i(f) & = & e_i(|\nu\cdot\Ga\ov{\Phi}_0-\ov{\Phi}_0|^2) \\
 & = & 2\,e_i\big(|\ov{\Phi}_0|^2+\mathrm{Re}\<\nu\cdot\Ga\ov{\Phi}_0,\ov{\Phi}_0\>\big)
 \end{eqnarray*}

Note that the spinor field $\Phi_0$ is parallel, so we can assume that $|\Phi_0|^2=1$ and since the trivialization is a fiberwise isometry, we have $|\ov{\Phi}_0|^2=1$ and then $e_i\big(|\ov{\Phi}_0|^2\big)=0$. Now using the compatibility of the Hermitian metric with the spinorial Levi-Civita connection and the properties~(\ref{poc}) of the chirality operator $\Ga$, we get:
\begin{eqnarray*}
e_i\big(\mathrm{Re}\<\nu\cdot\Ga\ov{\Phi}_0,\ov{\Phi}_0\>\big)=\mathrm{Re}\<\ov{\na}_{e_i}\nu\cdot\Ga\ov{\Phi}_0,\Phi_0\>+2\mathrm{Re}\<\nu\cdot\Ga(\ov{\na}_{e_i}\ov{\Phi}_0),\ov{\Phi}_0\>.
\end{eqnarray*}

However since the spinor field $\Phi_0$ is parallel, Formula~(\ref{levicivita}) leads to:
\begin{eqnarray*}
\ov{\na}_{e_i}\ov{\Phi}_0=\frac{1}{4}\sum_{j,k=1}^{n}\widetilde{\Ga}_{ij}^k e_j\cdot e_k\cdot\ov{\Phi}_0
\end{eqnarray*}

and thus we have $2\mathrm{Re}\<\nu\cdot\Ga(\ov{\na}_{e_i}\ov{\Phi}_0),\ov{\Phi}_0\>=\frac{1}{2}\sum_{j,k=1}^{n}\widetilde{\Ga}_{ij}^k\,\mathrm{Re}\<\nu\cdot\Ga (e_j\cdot e_k\cdot\ov{\Phi}_0),\ov{\Phi}_0\>$. We can now split this sum with respect to the tangent part and the normal part. Indeed, we write:
\begin{eqnarray*}
\sum_{j,k=1}^{n}\widetilde{\Ga}_{ij}^k e_j\cdot e_k\cdot\ov{\Phi}_0=\sum_{j,k=1}^{n-1}\widetilde{\Ga}_{ij}^k e_j\cdot e_k\cdot\ov{\Phi}_0-\sum_{j=1}^{n-1}\widetilde{\Ga}_{in}^j e_j\cdot \nu\cdot\ov{\Phi}_0+\sum_{j=1}^{n-1}\widetilde{\Ga}_{ij}^n e_j\cdot \nu\cdot\ov{\Phi}_0-\widetilde{\Ga}_{in}^n\ov{\Phi}_0
\end{eqnarray*}

and since in Fermi coordinates we have $\widetilde{\Ga}_{ij}^n=-\widetilde{\Ga}_{in}^j$ and $\widetilde{\Ga}^n_{in}=0$, we obtain:
\begin{eqnarray*}
\sum_{j,k=1}^{n}\widetilde{\Ga}_{ij}^k e_j\cdot e_k\cdot\ov{\Phi}_0=\sum_{j,k=1}^{n-1}\widetilde{\Ga}_{ij}^k e_j\cdot e_k\cdot\ov{\Phi}_0+2\sum_{j=1}^{n-1}\widetilde{\Ga}_{ij}^n e_j\cdot\nu\cdot\ov{\Phi}_0.
\end{eqnarray*}

This leads to:
\begin{eqnarray*}
\mathrm{Re}\<\nu\cdot\Ga(\ov{\na}_{e_i}\ov{\Phi}_0),\ov{\Phi}_0\> & = & \mathrm{Re}\<\ov{\na}_{e_i}\nu\cdot\Ga\ov{\Phi}_0,\ov{\Phi}_0\>+\frac{1}{2}\sum_{1\leq j\neq k\leq n-1}\widetilde{\Ga}_{ij}^k\,\mathrm{Re}\<\nu\cdot\Ga(e_j\cdot e_k\cdot\ov{\Phi}_0),\ov{\Phi}_0\>\\
& & -\frac{1}{2}\Big(\sum_{j=1}^{n-1}\widetilde{\Ga}_{ij}^j\Big)\,\mathrm{Re}\<\nu\cdot\Ga\ov{\Phi}_0,\ov{\Phi}_0\>+\sum_{j=1}^{n-1}\widetilde{\Ga}_{ij}^n\,\mathrm{Re}\<\nu\cdot\Ga (e_j\cdot\nu\cdot\ov{\Phi}_0),\ov{\Phi}_0\>
\end{eqnarray*}

and since $\widetilde{\Ga}_{ij}^j=0$ and $\widetilde{\Ga}_{ij}^n e_j=-\ov{\na}_{e_i}\nu$, we conclude that:
\begin{eqnarray*}
e_i(f) & = & \sum_{1\leq j\neq k\leq n-1}\widetilde{\Ga}_{ij}^k\,\mathrm{Re}\<\nu\cdot\Ga(e_j\cdot e_k\cdot\ov{\Phi}_0),\ov{\Phi}_0\>.
\end{eqnarray*}

An easy computation using the properties of the Hermitian metric and those of the chirality operator show that:
\begin{eqnarray*}
\<\nu\cdot\Ga(e_j\cdot e_k\cdot\ov{\Phi}_0),\ov{\Phi}_0\> & = & -\<\nu\cdot\Ga(e_j\cdot e_k\cdot\ov{\Phi}_0),\ov{\Phi}_0\>.
\end{eqnarray*}

Thus we have $\mathrm{Re}\<\nu\cdot\Ga(e_j\cdot e_k\cdot\ov{\Phi}_0),\ov{\Phi}_0\>=0$ and finally $e_i(f)=0$ for all $1\leq i\leq n-1$. Moreover since $f(q)=0$, the function $f$ vanishes identically on $\pa\V\cap\M$, i.e.: 
\begin{eqnarray*}
\nu\cdot\Ga\ov{\Phi}_{0|\pa\V\cap\M}=\ov{\Phi}_{0|\pa\V\cap\M}.
\end{eqnarray*}
 \hfill$\square$\\

We are now ready to prove the main theorem of this paper.\\
 
{\it Proof of Theorem~\ref{main2}:} Using Proposition~\ref{testspinor}, there exists a spinor field $\psi\in\Ga\big(\Si_{\xi}(\mathbb{R}^n_+)\big)$ satisfying:
\begin{eqnarray*}
\D_\xi\psi = \frac{n}{2} f\,\psi 
\end{eqnarray*}

where $f(x)=\frac{2}{1+r^2}$ and $r^2=x_1^2+...+x_{n-1}^2+t^2$, $|\psi|=f^{\frac{n-1}{2}}$ and $|\D_\xi\psi|= f^{\frac{n+1}{2}}$. Recall that this spinor field is given by:
\begin{eqnarray*}
\psi(x)=\frac{1}{\sqrt{2}}f^{\frac{n}{2}}(x)(1-x)\cdot\Phi_0(x)
\end{eqnarray*}

where $\Phi_0$ is a parallel spinor which is chosen to satisfy:
\begin{eqnarray}\label{pointboundary}
\nu\cdot\Ga\ov{\Phi}_0(q)=\ov{\Phi}_0(q),
\end{eqnarray}
 
at one point $q\in\V\cap\pa\M$. From this spinor field, we construct an adapted spinor on $\M$ for our problem. So let $\eta$ be a cut-off function given by $\eta=1$ on $B^+(q,\delta)$, $\eta=0$ on $\M\setminus B^+(q,2\delta)$ where $\delta$ is a positive number such that $\delta\leq 1$ and $B^+(q,z)$ is the half-ball of center $q$ and radius $z$ contained in the open set $\V$ defined in Section~\ref{spt}. Moreover, without loss of generality, we can assume that $\eta$ satisfies $|\ov{\na}{\eta}|\leq C r$, where $C$ is a positive real number. In the following, the symbol $C$ will stand for positive constants which can differ from one line to another. Let $\varepsilon>0$ be a small positive number, then we set:
\begin{eqnarray*}
\overline{\psi}_{\varepsilon}(x,t)=\eta\,\overline{\psi}\big(\frac{(x,t)}{\varepsilon}\big)\in\Gamma\big(\Sigma_g(\M)\big).
\end{eqnarray*}

Since the spinor field $\Phi_0$ satisfies (\ref{pointboundary}), Lemma~\ref{trivialisationbord} and the properties~(\ref{poc}) of the chirality operator lead to: 
\begin{eqnarray*}
 \nu\cdot\Ga\ov{\psi}_{\varepsilon|\pa\M}=\ov{\psi}_{\varepsilon|\pa\M}.
\end{eqnarray*}

Using Proposition~\ref{devdir} and Proposition~\ref{devdir1}, we have:
\begin{eqnarray*}
\D_g\ov{\psi}_{\varepsilon}\big(x,t\big) & = & \ov{\na}\eta\cdot\ov{\psi}\big(\frac{(x,t)}{\varepsilon}\big)+\eta\,\D_g\Big(\ov{\psi}\big(\frac{(x,t)}{\varepsilon}\big)\Big)\\ \\
& = & \ov{\na}\eta\cdot\ov{\psi}\big(\frac{(x,t)}{\varepsilon}\big)+\frac{\eta}{\varepsilon}\,\frac{n}{2}f\big(\frac{(x,t)}{\varepsilon}\big)\ov{\psi}\big(\frac{(x,t)}{\varepsilon}\big)+\frac{\eta}{\varepsilon}\sum_{i,j=1}^{n}\big(b_i^j-\delta_i^j\big)\ov{\pa_i\cdot\na_{\pa_j}\psi}\big(\frac{x}{\varepsilon}\big)
\\ \\ & & 
+\eta\,\mathrm{W}\cdot\ov{\psi}\big(\frac{(x,t)}{\varepsilon}\big)+\eta\,\mathrm{T}\cdot\ov{\psi}\big(\frac{(x,t)}{\varepsilon}\big)+\eta\,\wit{\nu}\cdot\mathrm{Z}\cdot\ov{\psi}\big(\frac{(x,t)}{\varepsilon}\big)\\ \\  & &
-\frac{n-1}{2}\eta\,\mathrm{H}_t\,\wit{\nu}\cdot\ov{\psi}\big(\frac{(x,t)}{\varepsilon}\big),
\end{eqnarray*}

with $|\mathrm{W}|=O(r^2)$, $|\mathrm{Z}|=O(r^2)$, $|\mathrm{T}|=O(r)$. Since the chiral bag invariant is independent of the metric chosen in the conformal class of $g$, we can assume that $g$ is such that $\R_g$ is of constant sign and $\H_g=0$. For this, one only needs to choose the conformal factor as being a solution of the eigenvalue problem~(\ref{pvpy}) and thus we obtain $|\mathrm{H}_t|=O(r)$. We now develop the third term of the preceding identity. In fact, we can easily check that:
\begin{eqnarray*}
\sum_{i,j=1}^{n}\big(b_i^j-\delta_i^j\big)\ov{\pa_i\cdot\na_{\pa_j}\psi}\big(\frac{(x,t)}{\varepsilon}\big) & = & \frac{1}{\sqrt{2}}f^{\frac{n}{2}}\sum_{i=1}^{n}\big(b_i^i-1\big)\ov{\Phi}_0-\frac{n}{\sqrt{2}}\,\ov{{\rm X}\cdot\psi}\big(\frac{(x,t)}{\varepsilon}\big),
\end{eqnarray*}

where ${\rm X}=f\,\sum_{i,j=1}^n\big(b_i^j-\delta_i^j\big)x_j\pa_i\in\Gamma(\TM)$. Using the development of $b^j_i$ given in (\ref{DEVBij}) yields to:
\begin{eqnarray}\label{bijdev}
\frac{1}{\sqrt{2}}f^{\frac{n}{2}}\sum_{i=1}^{n}\big(b_i^i-1\big)\ov{\Phi}_0 & = & \frac{1}{\sqrt{2}}f^{\frac{n}{2}}\big(O(r)\big)\ov{\Phi}_0
\end{eqnarray}

Then we develop:
\begin{eqnarray*}
|\D_g\overline{\psi}_{\varepsilon}|^2(x,t) & = & (I)+(II)+(III)+(IV)+(V)+(VI)+(VII)+(VIII)+(IX)+(X)\\
& & +(XI)+(XII)+(XIII)+(XIV)+(XV)+(VI)+(XVII)+(XVIII)\\
& & +(XIX)+(XX)+(XXI)+(XXII)+(XXIII)+(XXIV)+(XXV)\\
 & & +(XXVI)+(XXVII)+(XXVIII),
\end{eqnarray*}

where:
$$\begin{array}{lll}
(I) & = & |\ov{\na}\eta|^2|\ov{\psi}|^2\big(\frac{(x,t)}{\varepsilon}\big)\\ \\
(II) & = & \frac{\eta^2}{\varepsilon^2}\frac{n^2}{4}f^2\big(\frac{(x,t)}{\varepsilon}\big)|\ov{\psi}|^2\big(\frac{(x,t)}{\varepsilon}\big)\\ \\
(III) & = &
\frac{\eta^2}{\varepsilon^2}|\sum_{i,j=1}^{n}\big(b_i^j-\delta_i^j\big)\ov{\pa_i\cdot\na_{\pa_j}\psi}|^2\big(\frac{(x,t)}{\varepsilon}\big)\\ \\
(IV) & = & \eta^2|\mathrm{W}|^2|\ov{\psi}|^2\big(\frac{(x,t)}{\varepsilon}\big)\\ \\
(V) & = & \eta^{2}|\T|^2|\ov{\psi}|^2\big(\frac{(x,t)}{\varepsilon}\big)\\ \\
(VI) & = & \eta^{2}|\mathrm{Z}|^2|\ov{\psi}|^2\big(\frac{(x,t)}{\varepsilon}\big)\\ \\
(VII) & = & \eta^{2}h^2|\ov{\psi}|^2\big(\frac{(x,t)}{\varepsilon}\big)\\ \\
(VIII) & = &
 \frac{n\eta}{\varepsilon}f\big(\frac{(x,t)}{\varepsilon}\big)\mathrm{Re}\<\ov{\na}\eta\cdot\ov{\psi}\big(\frac{(x,t)}{\varepsilon}\big),\ov{\psi}\big(\frac{(x,t)}{\varepsilon}\big)\>\\ \\
(IX) & = & 
\frac{2\eta}{\varepsilon}\sum_{i,j=1}^{n}\big(b_i^j-\delta_i^j\big)\mathrm{Re}\<\ov{\na}\eta\cdot\ov{\psi}\big(\frac{(x,t)}{\varepsilon}\big),\ov{\pa_i\cdot\na_{\pa_j}\psi}\big(\frac{(x,t)}{\varepsilon}\big)\>\\ \\
\end{array}$$

$$\begin{array}{lll}
(X) & = & 2\eta\mathrm{Re}\<\ov{\na}\eta\cdot\ov{\psi}\big(\frac{(x,t)}{\varepsilon}\big),\mathrm{W}\cdot\ov{\psi}\big(\frac{(x,t)}{\varepsilon}\big)\>\\ \\
(XI) & = & 2\eta\mathrm{Re}\<\ov{\na}\eta\cdot\ov{\psi}\big(\frac{(x,t)}{\varepsilon}\big),\T\cdot\ov{\psi}\big(\frac{(x,t)}{\varepsilon}\big)\>\\ \\
(XII) & = & 2\eta\mathrm{Re}\<\ov{\na}\eta\cdot\ov{\psi}\big(\frac{(x,t)}{\varepsilon}\big),\wit{\nu}\cdot\mathrm{Z}\cdot\ov{\psi}\big(\frac{(x,t)}{\varepsilon}\big)\>\\ \\
(XIII) & = & 2 h\,\eta\mathrm{Re}\<\ov{\na}\eta\cdot\ov{\psi}\big(\frac{(x,t)}{\varepsilon}\big),\wit{\nu}\cdot\ov{\psi}\big(\frac{(x,t)}{\varepsilon}\big)\>\\ \\
(XIV) & = & \frac{n\eta}{\varepsilon^2}f\big(\frac{(x,t)}{\varepsilon}\big)\sum_{i,j=1}^n\big(b^j_i-\delta_i^j\big)\mathrm{Re}\<\ov{\psi}\big(\frac{(x,t)}{\varepsilon}\big),\ov{\pa_i\cdot\na_{\pa_j}\psi}\big(\frac{(x,t)}{\varepsilon}\big)\>\\ \\
(XV) & = & 
n\frac{\eta^2}{\varepsilon}f\big(\frac{(x,t)}{\varepsilon}\big)\mathrm{Re}\<\ov{\psi}\big(\frac{(x,t)}{\varepsilon}\big),\mathrm{W}\cdot\ov{\psi}\big(\frac{(x,t)}{\varepsilon}\big)\>\\ \\
(XVI) & = &  n\frac{\eta^2}{\varepsilon}f\big(\frac{(x,t)}{\varepsilon}\big)\mathrm{Re}\<\ov{\psi}\big(\frac{(x,t)}{\varepsilon}\big),\T\cdot\ov{\psi}\big(\frac{(x,t)}{\varepsilon}\big)\>\\ \\
(XVII) & = &  n\frac{\eta^2}{\varepsilon}f\big(\frac{(x,t)}{\varepsilon}\big)\mathrm{Re}\<\ov{\psi}\big(\frac{(x,t)}{\varepsilon}\big),\wit{\nu}\cdot\mathrm{Z}\cdot\ov{\psi}\big(\frac{(x,t)}{\varepsilon}\big)\>\\ \\
(XVIII) & = &  n\frac{\eta^2}{\varepsilon}f\big(\frac{(x,t)}{\varepsilon}\big)\mathrm{Re}\<\ov{\psi}\big(\frac{(x,t)}{\varepsilon}\big),h\,\wit{\nu}\cdot\ov{\psi}\big(\frac{(x,t)}{\varepsilon}\big)\>\\ \\
(XIX) & = & 
\frac{\eta^2}{\varepsilon}\sum_{i,j=1}^n\big(b^j_i-\delta_i^j\big)\mathrm{Re}\<\ov{\pa_i\cdot\na_{\pa_j}\psi}\big(\frac{(x,t)}{\varepsilon}\big),\mathrm{W}\cdot\ov{\psi}\big(\frac{(x,t)}{\varepsilon}\big)\>\\ \\

(XX) & = & 
\frac{\eta^2}{\varepsilon}\sum_{i,j=1}^n\big(b^j_i-\delta_i^j\big)\mathrm{Re}\<\ov{\pa_i\cdot\na_{\pa_j}\psi}\big(\frac{(x,t)}{\varepsilon}\big),\mathrm{T}\cdot\ov{\psi}\big(\frac{(x,t)}{\varepsilon}\big)\>\\ \\
(XXI) & = & 
\frac{\eta^2}{\varepsilon}\sum_{i,j=1}^n\big(b^j_i-\delta_i^j\big)\mathrm{Re}\<\ov{\pa_i\cdot\na_{\pa_j}\psi}\big(\frac{(x,t)}{\varepsilon}\big),\wit{\nu}\cdot\mathrm{Z}\cdot\ov{\psi}\big(\frac{(x,t)}{\varepsilon}\big)\>\\ \\
(XXII) & = & 
\frac{\eta^2}{\varepsilon}h\,\sum_{i,j=1}^n\big(b^j_i-\delta_i^j\big)\mathrm{Re}\<\ov{\pa_i\cdot\na_{\pa_j}\psi}\big(\frac{(x,t)}{\varepsilon}\big),\wit{\nu}\cdot\ov{\psi}\big(\frac{(x,t)}{\varepsilon}\big)\>\\ \\
(XXIII) & = & 2\eta^2\mathrm{Re}\<\mathrm{W}\cdot\ov{\psi}\big(\frac{(x,t)}{\varepsilon}\big),\T\cdot\ov{\psi}\big(\frac{(x,t)}{\varepsilon}\big)\>\\ \\
(XIV) & = & 2\eta^2\mathrm{Re}\<\mathrm{W}\cdot\ov{\psi}\big(\frac{(x,t)}{\varepsilon}\big),\wit{\nu}\cdot\mathrm{Z}\cdot\ov{\psi}\big(\frac{(x,t)}{\varepsilon}\big)\>\\ \\
(XV) & = & 2h\eta^2\mathrm{Re}\<\mathrm{W}\cdot\ov{\psi}\big(\frac{(x,t)}{\varepsilon}\big),\,\wit{\nu}\cdot\ov{\psi}\big(\frac{(x,t)}{\varepsilon}\big)\>\\ \\
(XVI) & = & 2\eta^2\mathrm{Re}\<\T\cdot\ov{\psi}\big(\frac{(x,t)}{\varepsilon}\big),\wit{\nu}\cdot\mathrm{Z}\cdot\ov{\psi}\big(\frac{(x,t)}{\varepsilon}\big)\>\\ \\
(XVII) & = & 2h\eta^2\mathrm{Re}\<\T\cdot\ov{\psi}\big(\frac{(x,t)}{\varepsilon}\big),\,\wit{\nu}\cdot\ov{\psi}\big(\frac{(x,t)}{\varepsilon}\big)\>\\ \\
(XVIII) & = & 2h\eta^2\mathrm{Re}\<\wit{\nu}\cdot\mathrm{Z}\cdot\ov{\psi}\big(\frac{(x,t)}{\varepsilon}\big),\,\wit{\nu}\cdot\ov{\psi}\big(\frac{(x,t)}{\varepsilon}\big)\>\\ \\
\end{array}$$

where $h=-\frac{n-1}{2}\H_t$. Since $\ov{\na}\eta$ and $\T$ are $1$-forms and $\mathrm{Z}$ is a $2$-form then:
$$\begin{array}{lll}
(VIII)=0,\qquad(XVI)=0,\qquad (XVIII)=0\qquad (XXVII)=0\qquad\text{and}\qquad (XVIII)=0.
\end{array}$$

Using this development, the properties of the fields $\mathrm{W}$, $\T$, $\U$ and $h$ given in Proposition \ref{devdir1} and since we assumed that $|\ov{\na}\eta|\leq Cr$ and $r\leq\delta\leq 1$, one can check that:
$$\begin{array}{ll}
\bullet & (I)+(IV)+(V)+(VI)+(VII)+(X)+(XI)+(XII)+(XIII)+(XXIII)+(XXIV)\\ \\
 & +(XXV)+(XXVI)\leq Cr^2f^{n-1}\big(\frac{(x,t)}{\varepsilon}\big)\\ \\
\bullet & (IX)+(XV)+(XVII)+(XIX)+(XX)+(XXI)+(XXII)\leq  \frac{C}{\varepsilon}r^2f^{n}\big(\frac{(x,t)}{\varepsilon}\big)+\frac{C}{\varepsilon}r^2f^{n+1}\big(\frac{(x,t)}{\varepsilon}\big)\\ \\
\bullet & (II)+(III)+(XIV)\leq\frac{n^2}{4\varepsilon^2}f^{n+1}\big(\frac{(x,t)}{\varepsilon}\big)+\frac{C}{\varepsilon^2}rf^{n+1}\big(\frac{(x,t)}{\varepsilon}\big)+\frac{C}{\varepsilon^2}r^2f^{n}\big(\frac{(x,t)}{\varepsilon}\big).
\end{array}$$

We can then write:
\begin{eqnarray*}
0\leq |\D_g\overline{\psi}_{\varepsilon}|^2(x,t)& \leq & \frac{n^2}{4\varepsilon^2}f^{n+1}\big(\frac{(x,t)}{\varepsilon}\big)+\frac{C}{\varepsilon^2}rf^{n+1}\big(\frac{(x,t)}{\varepsilon}\big)+\frac{C}{\varepsilon^2}r^2f^{n}\big(\frac{(x,t)}{\varepsilon}\big)\\ \\ 
& & +\frac{C}{\varepsilon}r^2f^{n}\big(\frac{(x,t)}{\varepsilon}\big) +\frac{C}{\varepsilon}r^2f^{n+1}\big(\frac{(x,t)}{\varepsilon}\big)+Cr^2f^{n-1}\big(\frac{(x,t)}{\varepsilon}\big)\\ \\
& \leq & \frac{n^2}{4\varepsilon^2}f^{n+1}\big(\frac{(x,t)}{\varepsilon}\big)\big[1+\Lambda\big],
\end{eqnarray*}

where $\Lambda=C r+C r^2f^{-1}\big(\frac{(x,t)}{\varepsilon}\big)+C\varepsilon r^2+C r^2\varepsilon f^{-1}\big(\frac{(x,t)}{\varepsilon}\big)+C r^2\varepsilon^2f^{-2}\big(\frac{(x,t)}{\varepsilon}\big)$. Note that since $|\D_g\overline{\psi}_{\varepsilon}|^2\geq 0$, then $\Lambda\geq -1$. On the other hand, if $x\geq -1$ we have:
\begin{eqnarray*}
(1+x)^{\frac{n}{n+1}}\leq 1+\frac{n}{n+1}x,
\end{eqnarray*}

so we get:
\begin{eqnarray*}
|\D_g\overline{\psi}_{\varepsilon}|^{\frac{2n}{n+1}}(x,t) & \leq & \Big(\frac{n^2}{4\varepsilon^2}f^{n+1}\big(\frac{(x,t)}{\varepsilon}\big)\Big)^{\frac{n}{n+1}}\big[1+\Lambda\big]^{\frac{n}{n+1}}\\
& \leq & \Big(\frac{n}{2\varepsilon}\Big)^{\frac{2n}{n+1}}f^n\big(\frac{(x,t)}{\varepsilon}\big)\big[1+\frac{n}{n+1}\Lambda\big],\end{eqnarray*}

i.e.:
\begin{eqnarray*}
|\D_g\overline{\psi}_{\varepsilon}|^{\frac{2n}{n+1}}(x,t) & \leq & \Big(\frac{n}{2}\Big)^{\frac{2n}{n+1}}\varepsilon^{-\frac{2n}{n+1}}\Big[f^n\big(\frac{(x,t)}{\varepsilon}\big)+C\frac{n}{n+1}rf^n\big(\frac{(x,t)}{\varepsilon}\big)\\ \\ & & +C\frac{n}{n+1}r^2f^{n-1}\big(\frac{(x,t)}{\varepsilon}\big) +C\frac{n}{n+1}\varepsilon r^2 f^{n}\big(\frac{(x,t)}{\varepsilon}\big)\\ \\
& & +C\frac{n}{n+1}\varepsilon r^2f^{n-1}\big(\frac{(x,t)}{\varepsilon}\big)+C\frac{n}{n+1}\varepsilon^2 r^2f^{n-2}\big(\frac{(x,t)}{\varepsilon}\big)\Big].
\end{eqnarray*}

Integrating the last inequality leads to:
\begin{eqnarray*}
\int_{\M}|\D_g\overline{\psi}_{\varepsilon}|^{\frac{2n}{n+1}}dv(g) & \leq & \varepsilon^{-\frac{2n}{n+1}}\Big[{\bf A}+{\bf B}+{\bf C}+{\bf D}+{\bf E}+{\bf F}\Big],
\end{eqnarray*}

where:
\begin{eqnarray*}
{\bf A} & = & \int_{\mathrm{B}_q^+(2\delta)}f^n\big(\frac{(x,t)}{\varepsilon}\big)dv(g) \\ \\
{\bf B} & = & C\frac{n}{n+1}\int_{\mathrm{B}_q^+(2\delta)} rf^n\big(\frac{(x,t)}{\varepsilon}\big)dv(g) \\ \\
{\bf C} & = & C\frac{n}{n+1}\int_{\mathrm{B}_q^+(2\delta)}r^2f^{n-1}\big(\frac{(x,t)}{\varepsilon}\big)dv(g) \\ \\
{\bf D} & = & C\frac{n}{n+1}\varepsilon \int_{\mathrm{B}_q^+(2\delta)}r^2 f^{n}\big(\frac{(x,t)}{\varepsilon}\big)dv(g)\\ \\
{\bf E} & = & C\frac{n}{n+1}\varepsilon \int_{B_q^+(2\delta)}r^2f^{n-1}\big(\frac{(x,t)}{\varepsilon}\big)dv(g)\\ \\
{\bf F} & = & C\frac{n}{n+1}\varepsilon^2 \int_{\mathrm{B}_q^+(2\delta)}r^2f^{n-2}\big(\frac{(x,t)}{\varepsilon}\big)dv(g).
\end{eqnarray*}

Now we have to estimate these terms. Let's start with ${\bf A}$. Note that since the function $f$ is radial, spherical coordinates lead to:
\begin{eqnarray*}
{\bf A} & = & \frac{\omega_{n-1}}{2}\int_{\mathrm{B}_q^+(2\delta)}r^{n-1}f^n\big(\frac{(x,t)}{\varepsilon}\big)\mathcal{S}(r)dr,
\end{eqnarray*}

where:
\begin{eqnarray*}
\mathcal{S}(r)=\frac{2}{\omega_{n-1}}\int_{\mathbb{S}^{n-1}_+}\sqrt{\mathrm{det}(g_{rx})}ds(x).
\end{eqnarray*}

However with the help of the development of the volume form of $(\M,g)$ (see \cite{esc}), we have:
\begin{eqnarray*}
\mathcal{S}(r)\leq 1+Cr,
\end{eqnarray*}

and thus:
\begin{eqnarray*}
{\bf A} & \leq & \frac{\omega_{n-1}}{2}\Big[\int_{\mathrm{B}_q^+(2\delta)}r^{n-1}f^n\big(\frac{(x,t)}{\varepsilon}\big)dr+C\int_{\mathrm{B}_q^+(2\delta)}r^{n}f^n\big(\frac{(x,t)}{\varepsilon}\big)dr\Big].
\end{eqnarray*}
 
A simple change of variables gives:
\begin{eqnarray*}
{\bf A} & \leq & \frac{\omega_{n-1}}{2}\varepsilon^n \Big[\int_{0}^{\frac{2\delta}{\varepsilon}}r^{n-1}f^n(r)dr+C\varepsilon\int_{0}^{\frac{2\delta}{\varepsilon}}r^{n}f^n(r)dr\Big].
\end{eqnarray*}

Some calculations show that if $n\geq 1$: 
\begin{eqnarray*}
\int_{\frac{2\delta}{\varepsilon}}^{+\infty}r^{n-1}f^n(r)dr & = & o(1),\\ \\
\int_{0}^{\frac{2\delta}{\varepsilon}}r^n f^n(r)dr & = & o(\mathrm{ln}\varepsilon).
\end{eqnarray*}

We finally get:
\begin{eqnarray*}
{\bf A} & \leq & \frac{\omega_{n-1}}{2}\varepsilon^n \Big[\int_{0}^{+\infty}r^{n-1}f^n(r)dr+o(1)\Big].
\end{eqnarray*}

We now give an upper bound for ${\bf B}$. In the same way, we compute:
\begin{eqnarray*}
{\bf B} & = &  C\int_{\mathrm{B}_q^+(2\delta)} rf^n\big(\frac{(x,t)}{\varepsilon}\big)dv(g)\\
 & \leq & C\int_{\mathrm{B}_q^+(2\delta)} rf^n\big(\frac{(x,t)}{\varepsilon}\big)dx+C\int_{\mathrm{B}_q^+(2\delta)} r^2f^n\big(\frac{(x,t)}{\varepsilon}\big)dx
\end{eqnarray*}

i.e.:
\begin{eqnarray*}
{\bf B} & \leq & C\varepsilon^{n+1}\int^{\frac{2\delta}{\varepsilon}}_0 r^n f^n(r)dr+C\varepsilon^{n+2}\int^{\frac{2\delta}{\varepsilon}}_0 r^{n+1}f^n(r)dr.
\end{eqnarray*}

We can easily check that if $n\geq 1$ then ${\bf B}=o(\varepsilon^n)$. For the quantity ${\bf C}$, we write:
\begin{eqnarray*}
{\bf C} & = & C\int_{\mathrm{B}_q^+(2\delta)}r^2f^{n-1}\big(\frac{(x,t)}{\varepsilon}\big)dv(g) \\
 & \leq & C\int_{\mathrm{B}_q^+(2\delta)}r^2f^{n-1}\big(\frac{(x,t)}{\varepsilon}\big)dx+C\int_{\mathrm{B}_q^+(2\delta)}r^3f^{n-1}\big(\frac{(x,t)}{\varepsilon}\big)dx
\end{eqnarray*}

which after a change of variable gives:
\begin{eqnarray*}
{\bf C} & \leq & C\varepsilon^{n+2}\int^{\frac{2\delta}{\varepsilon}}_0r^{n+1}f^{n-1}(r)dr+C\varepsilon^{n+3}\int^{\frac{2\delta}{\varepsilon}}_0r^{n+2}f^{n-1}(r)dr.
\end{eqnarray*}

We can thus conclude that if $n\geq 3$, ${\bf C}=o(\varepsilon^n)$. The term ${\bf D}$ satisfies:
\begin{eqnarray*}
{\bf D} & = & C\varepsilon \int_{\mathrm{B}_q^+(2\delta)}r^2 f^{n}\big(\frac{(x,t)}{\varepsilon}\big)dv(g)\\ 
& \leq & C\varepsilon \int_{\mathrm{B}_q^+(2\delta)}r^2 f^{n}\big(\frac{(x,t)}{\varepsilon}\big)dx+C\varepsilon \int_{\mathrm{B}_q^+(2\delta)}r^3 f^{n}\big(\frac{(x,t)}{\varepsilon}\big)dx
\end{eqnarray*}

and then:
\begin{eqnarray*}
{\bf D} & \leq & C\varepsilon^{n+2} \int^{\frac{2\delta}{\varepsilon}}_0r^{n+1} f^{n}(r)dr+C\varepsilon^{n+3} \int^{\frac{2\delta}{\varepsilon}}_0r^{n+2} f^{n}(r)dr.
\end{eqnarray*}

So we have shown that if $n\geq 1$, ${\bf D}=o(\varepsilon^n)$. In the same way, we compute that if $n\geq 3$, then ${\bf E}$ and ${\bf F}$ satisfy ${\bf E}=o(\varepsilon^n)$ and ${\bf F}=o(\varepsilon^n)$. Finally we get that for $n\geq 3$:
\begin{eqnarray*}
\int_{\M}|\D_g\overline{\psi}_{\varepsilon}|^{\frac{2n}{n+1}}dv(g) & \leq & \Big(\frac{n}{2}\Big)^{\frac{2n}{n+1}}\frac{\omega_{n-1}}{2}\varepsilon^{\frac{n(n-1)}{n+1}}\Big[\int_{0}^{\infty}r^{n-1}f^n(r)dr+o(1)\Big],
\end{eqnarray*}

and so:
\begin{eqnarray}
\Big(\int_{\M}|\D_g\overline{\psi}_{\varepsilon}|^{\frac{2n}{n+1}}dv(g)\Big)^{\frac{n+1}{n}}\leq \frac{n^2}{4}\Big(\frac{\omega_{n-1}}{2}\Big)^{\frac{n+1}{n}}\;\mathrm{I}^{\frac{n+1}{n}}\varepsilon^{n-1}\big(1+o(1)\big),
\end{eqnarray}

where $\mathrm{I}:=\int_{0}^{\infty}r^{n-1}f^n(r)dr$. We are now going to estimate the denominator of the variational characterization of $\lambda_{\min}(\M,\pa\M)$. Indeed we have:
\begin{eqnarray*}
\mathrm{Re}\<\D_g\overline{\psi}_{\varepsilon},\overline{\psi}_{\varepsilon}\> & = & (I')+(II')+(III')+(IV')+(V')+(VI')+(VII'),
\end{eqnarray*} 

where:
\begin{eqnarray*}
(I') & = & \eta\mathrm{Re}\<\ov{\na}\eta\cdot\ov{\psi}\big(\frac{(x,t)}{\varepsilon}\big),\ov{\psi}\big(\frac{(x,t)}{\varepsilon}\big)\>\\ \\
(II') & = & \frac{n}{2\varepsilon}\eta^2\mathrm{Re}\<f\big(\frac{(x,t)}{\varepsilon}\big)\ov{\psi}\big(\frac{(x,t)}{\varepsilon}\big),\ov{\psi}\big(\frac{(x,t)}{\varepsilon}\big)\>\\ \\
(III') & = & 
\frac{1}{\varepsilon}\sum_{i,j=1}^n\eta^2\big(b^j_i-\delta_i^j\big)\mathrm{Re}\<\ov{\pa_i\cdot\na_{\pa_j}\psi}\big(\frac{(x,t)}{\varepsilon}\big),\ov{\psi}\big(\frac{(x,t)}{\varepsilon}\big)\>
\\ \\
(IV') & = & \eta^2\mathrm{Re}\<\mathrm{W}\cdot\ov{\psi}\big(\frac{(x,t)}{\varepsilon}\big),\ov{\psi}\big(\frac{(x,t)}{\varepsilon}\big)\>\\ \\
(V') & = & \eta^2\mathrm{Re}\<\wit{\nu}\cdot\mathrm{Z}\cdot\ov{\psi}\big(\frac{(x,t)}{\varepsilon}\big),\ov{\psi}\big(\frac{(x,t)}{\varepsilon}\big)\> \\ \\
(VI') & = & \eta^2\mathrm{Re}\<\T\cdot\ov{\psi}\big(\frac{(x,t)}{\varepsilon}\big),\ov{\psi}\big(\frac{(x,t)}{\varepsilon}\big)\>\\ \\
(VII') & = & \eta^2\mathrm{Re}\< h\wit{\nu}\cdot\ov{\psi}\big(\frac{(x,t)}{\varepsilon}\big),\ov{\psi}\big(\frac{(x,t)}{\varepsilon}\big)\>.
\end{eqnarray*} 

We then easily check that since:
\begin{eqnarray*}
(I')=0,\quad (VI')=0\quad\textrm{and}\quad (VII')=0
\end{eqnarray*}

we have:
\begin{eqnarray*}
\Big|\int_{\M}\mathrm{Re}\<\D_g\overline{\psi}_{\varepsilon},\overline{\psi}_{\varepsilon}\>dv(g)\Big|= {\bf A'}+{\bf B'}+{\bf C'}+{\bf D'},
\end{eqnarray*}

where we let:
\begin{eqnarray*}
{\bf A'} & = & \frac{n}{2\varepsilon}\int_{\mathrm{B}_q^+(2\delta)}f\big(\frac{(x,t)}{\varepsilon}\big)|\ov{\psi}|^2\big(\frac{(x,t)}{\varepsilon}\big)dv(g) =  \frac{n}{2\varepsilon}\int_{\mathrm{B}_q^+(2\delta)}f^n\big(\frac{(x,t)}{\varepsilon}\big)dv(g)\\ \\
{\bf B'} & = & \frac{C}{\varepsilon}\sum_{i,j=1}^n\int_{\mathrm{B}_q^{+}(2\delta)}\big(b^j_i-\delta_i^j\big)\mathrm{Re}\<\ov{\pa_i\cdot\na_{\pa_j}\psi}\big(\frac{(x,t)}{\varepsilon}\big),\ov{\psi}\big(\frac{(x,t)}{\varepsilon}\big)\>dv(g)\\
{\bf C'} & = &  \int_{\mathrm{B}_q^+(2\delta)}\mathrm{Re}\<\mathrm{W}\cdot\ov{\psi}\big(\frac{(x,t)}{\varepsilon}\big),\ov{\psi}\big(\frac{(x,t)}{\varepsilon}\big)\>dv(g)\\ \\
{\bf D'} & = & \int_{\mathrm{B}_q^+(2\delta)}\mathrm{Re}\<\wit{\nu}\cdot\mathrm{Z}\cdot\ov{\psi}\big(\frac{(x,t)}{\varepsilon}\big),\ov{\psi}\big(\frac{(x,t)}{\varepsilon}\big)\>dv(g).
\end{eqnarray*}

However using Proposition \ref{devdir1}, we have $|\mathrm{W}|=O(r^2)$ and $|\mathrm{Z}|=O(r^2)$, thus:
\begin{eqnarray*}
{\bf C'}+{\bf D'} & = & C\int_{\mathrm{B}_q^+(2\delta)}r^2f^{n-1}\big(\frac{(x,t)}{\varepsilon}\big)dv(g),
\end{eqnarray*}

so:
$$\begin{array}{ll}
\int_{\M}\mathrm{Re}\<\D_g\overline{\psi}_{\varepsilon},\overline{\psi}_{\varepsilon}\>& dv(g)   =  \frac{n}{2\varepsilon}\int_{\mathrm{B}_q^+(2\delta)}f^n\big(\frac{(x,t)}{\varepsilon}\big)dv(g)+C\int_{\mathrm{B}_q^+(2\delta)}r^2f^{n-1}\big(\frac{(x,t)}{\varepsilon}\big)dv(g)\\ \\
&  +\frac{C}{\varepsilon}\sum_{i,j=1}^n\int_{\mathrm{B}^{+}_q(2\delta)}\big(b^j_i-\delta_i^j\big)\mathrm{Re}\<\ov{\pa_i\cdot\na_{\pa_j}\psi}\big(\frac{(x,t)}{\varepsilon}\big),\ov{\psi}\big(\frac{(x,t)}{\varepsilon}\big)\>dv(g).
\end{array}$$

As done for the numerator, we compute:
\begin{eqnarray*}
{\bf A'} & = & \frac{n}{2\varepsilon}\int_{\mathrm{B}_q^+(2\delta)}f^n\big(\frac{(x,t)}{\varepsilon}\big)dv(g)\\
 & = & \frac{n}{2\varepsilon}\int_{\mathrm{B}_q^+(2\delta)}f^n\big(\frac{(x,t)}{\varepsilon}\big)dx+\frac{C}{\varepsilon}\int_{\mathrm{B}_q^+(2\delta)}rf^n\big(\frac{(x,t)}{\varepsilon}\big)dx
\end{eqnarray*}

and for $n\geq 1$, we conclude:
\begin{eqnarray*}
{\bf A'} & = & \frac{n}{4}\,\omega_{n-1}\mathrm{I}\,\varepsilon^{n-1} \big(1+o(1)\big).
\end{eqnarray*}

Using the estimate of the term ${\bf B}$, we check that:
\begin{eqnarray*}
{\bf B'} = \frac{C}{\varepsilon}\int_{\mathrm{B}^{+}_q(2\delta)}rf^{n}\big(\frac{(x,t)}{\varepsilon}\big)dv(g)
\end{eqnarray*}
 
and so ${\bf B'}=o(\varepsilon^{n-1})$ if $n\geq3$. An similar calculation shows that ${\bf C'}+{\bf D'}=o(\varepsilon^{n-1})$ if $n\geq 3$. Finally we get:
\begin{eqnarray*}
\Big|\int_{\M}\mathrm{Re}\<\D_g\overline{\psi}_{\varepsilon},\overline{\psi}_{\varepsilon}\>dv(g)\Big| & = & \frac{n}{2}\,\frac{\omega_{n-1}}{2}\mathrm{I}\,\varepsilon^{n-1} \big(1+o(1)\big).
\end{eqnarray*}

Now using the variational characterization of $\lambda_{\min}(\M,\pa\M)$ given in Proposition~\ref{vcbhl}, we obtain:
\begin{eqnarray*}
\lambda_{\min}(\M,\pa\M)\leq\frac{\big(\int_{\M}|\D_g\ov{\psi}_{\varepsilon}|^{\frac{2n}{n+1}}dv(g)\big)^{\frac{n+1}{n}}}{\big|\int_{\M}\mathrm{Re}\<\D_g\ov{\psi}_{\varepsilon},\ov{\psi}_{\varepsilon}\>dv(g)\big|}.
\end{eqnarray*}

The estimate of this functional allows to write:
\begin{eqnarray*}
\frac{\big(\int_{\M}|\D_g\ov{\psi}_{\varepsilon}|^{\frac{2n}{n+1}}dv(g)\big)^{\frac{n+1}{n}}}{\big|\int_{\M}\mathrm{Re}\<\D_g\ov{\psi}_{\varepsilon},\ov{\psi}_{\varepsilon}\>dv(g)\big|}\leq \frac{\frac{n^2}{4}\Big(\frac{\omega_{n-1}}{2}\Big)^{\frac{n+1}{n}}\mathrm{I}^{\frac{n+1}{n}}\varepsilon^{n-1}}{\frac{n}{2}\,\frac{\omega_{n-1}}{2}\,\mathrm{I}\,\varepsilon^{n-1}}\Big(1+o(1)\Big),
\end{eqnarray*}

which gives:
\begin{eqnarray*}
\lambda_{\min}(\M,\pa\M)\leq \frac{\big(\int_{\M}|\D_g\ov{\psi}_{\varepsilon}|^{\frac{2n}{n+1}}dv(g)\big)^{\frac{n+1}{n}}}{\big|\int_{\M}\mathrm{Re}\<\D_g\ov{\psi}_{\varepsilon},\ov{\psi}_{\varepsilon}\>dv(g)\big|}\leq \frac{n}{2}\Big(\frac{\omega_{n-1}}{2}\Big)^{\frac{1}{n}}\,\mathrm{I}^{\frac{1}{n}}\Big(1+o(1)\Big).
\end{eqnarray*}

However since $\omega_{n-1}\mathrm{I}=\omega_n$, we get:
\begin{eqnarray*}
\lambda_{\min}(\M,\pa\M)\leq\frac{n}{2}\Big(\frac{\omega_n}{2}\Big)^\frac{1}{n}\Big(1+o(1)\Big)=\lambda_{\min}(\hs,\pa\hs)\Big(1+o(1)\Big),
\end{eqnarray*}

and we can thus conclude that $\lambda_{\min}(\M,\pa\M)\leq\lambda_{\min}(\hs,\pa\hs)$.
\hfill$\square$\\

\begin{remark}
This result gives in particular a spinorial proof of Escobar's result given by (\ref{pablo}). Indeed, using the Hijazi inequality~(\ref{hijbord}), we note that:
\begin{eqnarray*}
\frac{n}{4(n-1)}\mu(\M,\pa\M)\leq\lambda_{\min}(\M,\pa\M)^2\leq\lambda_{\min}(\hs,\pa\hs)^2=\frac{n^2}{4}\Big(\frac{\omega_{n}}{2}\Big)^{\frac{2}{n}}
\end{eqnarray*}

and thus
\begin{eqnarray*}
\mu(\M,\pa\M)\leq n(n-1)\Big(\frac{\omega_{n}}{2}\Big)^{\frac{2}{n}}=\mu(\hs,\pa\hs).
\end{eqnarray*}
\end{remark}


\section{The case of surfaces with boundary}


We can show that Inequality~(\ref{largeb1}) still holds if $\M$ is a compact Riemannian surface with connected boundary. However, the test spinor needs a slight modification. We follow the argument given in \cite{gh} (see also \cite{amm3}). In fact, we show:
\begin{theorem}\label{insur}
Let $(\M^2,g)$ be a compact Riemannian surface with connected boundary. Then the chiral bag invariant satisfies:
\begin{eqnarray}\label{insur1}
\lambda_{\min}(\M,\pa\M)\leq \lambda_{\min}(\mathbb{S}^2_+,\pa\mathbb{S}^2_+)=\sqrt{2\pi}.
\end{eqnarray}
\end{theorem}

{\it Sketch of proof:} First note that since $\M$ is a surface it is always spin and it is equipped with a chirality operator given by the complex volume element of the spinor bundle. Moreover we can assume that $g$ is locally conformally flat and since the boundary is connected it is umbilic. Now let $0<\varepsilon\leq\alpha\leq\delta$ be real numbers and consider the positive function defined by:
$$f_{\varepsilon}(x)=
\left\lbrace
\begin{array}{ll} 
\frac{2\varepsilon^2}{\varepsilon^2+r^2} & \quad\rm{if}\;\;r\leq \alpha\\
\frac{2\varepsilon^2}{\varepsilon^2+\alpha^2} & \quad\rm{if}\;\;r\geq\alpha
\end{array}
\right.$$

where $r=d(q,x)$ and $q\in\pa\M$. Using Proposition~\ref{testspinor}, there exists a spinor field $\psi\in\Gamma\big(\Sigma_\xi(\mathbb{R}^2_+)\big)$ which satisfies:
\begin{eqnarray*}
\D_\xi\psi= f\psi 
\end{eqnarray*}

where $f(r)=\frac{2}{1+r^2}\in\mathrm{C}^{\infty}(\M)$ and the spinor $\psi$ can be written:
\begin{eqnarray*}
\psi=\frac{1}{\sqrt{2}}f\,(1-x)\cdot\Phi_0,
\end{eqnarray*} 

where $\Phi_0$ is a parallel spinor field such that $\mathbb{B}^-_g(\ov{\Phi}_0)(q)=0$. Now let $\eta$ a smooth function on $\M$ such that: 
$$\eta(x)=
\left\lbrace
\begin{array}{lll}
1 & \rm{on} & \mathrm{B}_q^+(\delta)\\
0 & \rm{on} & \M\setminus\mathrm{B}_q^+(2\delta) 
\end{array}
\right.
\qquad\rm{and}\qquad |\ov{\nabla}\eta|\leq\frac{1}{\delta}$$

with $\mathrm{B}_q^+(2\delta)\subset\V$ and $\V$ is an open flat subset of $\M$ defined in Section~\ref{spt}. We then consider the test spinor given by $\psi_{\varepsilon}(x)=\eta(x)\psi(\frac{x}{\varepsilon})$ (which by construction satisfies $\mathbb{B}^-_g(\ov{\psi}_{\varepsilon\,|\pa\M})=0$) and so using Corollary~\ref{pvd} we have:
\begin{eqnarray*}
\lambda_{1}(g_{\varepsilon})\leq \frac{\int_{\M}|\D_\varepsilon\ov{\psi}_{\varepsilon}|^2 f_{\varepsilon}^{-1}dv(g_\varepsilon)}{\big|\int_{\M}\mathrm{Re}\<\D_\varepsilon\ov{\psi}_{\varepsilon},\ov{\psi}_{\varepsilon}\>dv(g_\varepsilon)\big|}
\end{eqnarray*}

where $\lambda_{1}(g_{\varepsilon})$ is the first eigenvalue of the Dirac operator $\D_\varepsilon$ under the chiral bag boundary condition in the metric $g_{\varepsilon}=f_{\varepsilon}^2 g\in[g]$. We can then estimate this ratio and we obtain:
\begin{eqnarray*}
\lambda_{1}(g_{\varepsilon})\leq\frac{1}{\varepsilon}+o(1).
\end{eqnarray*}

Now we compute the volume of the surface $\M$ equipped with the metric $g_{\varepsilon}$; we clearly have:
\begin{eqnarray*}
\mathrm{Vol}(\M,g_{\varepsilon}) = \int_{\M}f_{\varepsilon}^2 dv(g) = \varepsilon^2\omega_1\big(1+o(\varepsilon)\big).
\end{eqnarray*}

Combining these estimations leads to: 
\begin{eqnarray*}
\lambda_{\min}(\M,\pa\M)\leq\lambda_{1}(g_{\varepsilon})\mathrm{Vol}(\M,g_{\varepsilon})^{\frac{1}{2}}=\omega_1^{\frac{1}{2}}\big(1+o(1)\big)=\sqrt{2\pi}\big(1+o(1)\big)
\end{eqnarray*}

and then Inequality~(\ref{insur1}) follows directly since $\lambda_{\min}(\mathbb{S}^2_+,\pa\mathbb{S}^2_+)=\sqrt{2\pi}$.
\hfill$\square$


\bibliographystyle{amsalpha}     
\bibliography{bibthese1}


\end{document}